\documentclass[twoside,a4paper,12pt]{book}

\usepackage{graphicx}
\usepackage{caption}
\usepackage{subcaption}

\graphicspath{Figures/}
\usepackage{comment}
\usepackage[absolute,overlay]{textpos}
\usepackage{graphicx,rotating,booktabs}
\usepackage{cite}
\usepackage{blindtext}
\usepackage{floatflt}
\usepackage{tikz}
\usepackage{setspace}
\usepackage{ragged2e}
\usepackage[hidelinks]{hyperref}
\usepackage{caption}
\usepackage{dirtytalk}
\usepackage{fancyhdr}
\usepackage[toc,page]{appendix}
\usepackage{pifont}
\usepackage[square, numbers, comma, sort&compress]{natbib}
\usepackage{verbatim}
\usepackage{./Extra_files/vector}
\usepackage{amssymb}
\usepackage{amsmath}
\usepackage{hyperref}
\usepackage{tocloft}
\usepackage{pdflscape}
\usepackage{color}
\usepackage[export]{adjustbox}
\usepackage{enumitem} 
\usepackage{multirow}
\usepackage{wrapfig} 
\usepackage{listings}

\numberwithin{equation}{section}
\numberwithin{equation}{section}

\newcommand{\blue}[1]{{\color{blue} #1}}

\newtheorem{theorem}{\bf \blue{Theorem}}[section]

\fancypagestyle{mystyle}
{
\fancyhf{}
\fancyhead[L]{\fontsize{8}{10} \selectfont \bfseries\rightmark}
\fancyhead[R]{\bfseries\thepage}

\addtolength{\headheight}{1pt}
\fancyfoot[L]{\bfseries \it{\bfseries{}}}
}

\begin{document}
\setstretch{.9}  
\pagenumbering{gobble}

\textbf{\Large{\begin{center} A Mathematical Modeling Study of COVID-19 With Reference to
Immigration from Urban to Rural Population \end{center}}}

\vspace{2cm}

{
D. K. K. Vamsi$^{a,1 }$, C. Bishal Chhetri$^{a}$, D. Bhanu prakash $^{a},$  Seshasainath Ch. $^{a},$  D. Surabhi Pandey $^{b,*}$  \\\\ \medskip 

$^{a}$Department of Mathematics and Computer Science, Sri Sathya Sai Institute of Higher Learning -  SSSIHL, India  \\ 
  
$^{b}$ Indian Institute of Public Health, Delhi  \\\\ 

dkkvamsi@sssihl.edu.in,  bishalchhetri@sssihl.edu.in, dbhanuprakash@sssihl.edu.in, seshasainath.ch@gmail.com, surabhi.pandey@phfi.org \\
{\small $^{1}$  First Author},
	{ \small $^{*}$ Corresponding Author}
}
\newpage
\begin{center}
{
\Large\bf\underline{Abstract}}
\end{center}

\textbf{\fontsize{13}{20}
{\flushleft{
In this study, we have formulated and analyzed a non-linear compartmental model (SEIR) for the dynamics of COVID-19 with reference to immigration from urban to rural population in Indian scenario. We have captured the effect of the immigration as two separate factors contributing in the rural compartments of the model.   We have first established the positivity of the solution and the boundedness of the solution followed by the existence and uniqueness of the solution for this multi compartment model. We later went on to find out the equibria of the system and derived the reproduction number. Further we numerically depicted the local and global stability of the equilibria. Later we have done sensitivity analysis of the model parameters and identified the sensitive parameters of the system. The sensitivity analysis is followed up with the two parameter  heat plots dealing with the sensitive parameters of the system. These heat plots gives us the  parameter regions in which the system is stable. Finally comparative and effectiveness studies were done with reference to the control interventions such as Vaccination, Antiviral drugs, Immunotheraphy. 
}} 
\\ \\
}

\newpage\null\thispagestyle{empty}

\fancyhead{}	
\begingroup
\let\cleardoublepage\clearpage
\tableofcontents
\endgroup

\oddsidemargin 1.36cm
\evensidemargin -0.0cm
\newpage\null\thispagestyle{empty}
\listoffigures
\newpage\null\thispagestyle{empty}

\listoftables

\phantomsection	
\addtocontents{lof}
{
\protect\setlength
{\protect\cftbeforefigskip}{7pt}
}

\addtocontents{toc}{\protect\hypertarget{toc}{}}
\let\cleardoublepage\clearpage

\addcontentsline{toc}{chapter}{\listfigurename}
\addcontentsline{toc}{chapter}{\listtablename}

\newcommand*\cleartoleftpage
{
\clearpage\ifodd\value{page}\hbox{}\newpage\fi
}

\cleartoleftpage
\pagestyle{mystyle}			

\pagenumbering{arabic}		

\chapter[Introduction and Motivation   ] {\hyperlink{toc}{Introduction and Motivation   }}
\setlength{\headheight}{14.49998pt}
\thispagestyle{empty}

\section{Brief Overview of Epidemic Modeling}

It is a well recognized fact that epidemic outbreaks in a community affect the lives of thousands of people. Due to high mortality and morbidity rates and the various disease-related costs such as expenditure on health care, diagnosis etc., the economy of the community is heavily disturbed. Disease like influenza, flu, SARS etc. have majorly contributed at the global level for these causes. Thus the prevention and control of any infectious disease has become utmost essential. 
\vspace{.2cm} \\

Mathematical Models have been used extensively to control, predict and formulate policies so as to eradicate the epidemic outbreaks and to study the disease burden {\cite{BIF, GLB}}.  
\vspace{.2cm} \\

One of the factors in reducing an infectious disease burden is the rate at which people get recovered which in turn depends on the number of individuals in the infected class. The relation between various compartments and the population level at any given point of time is the basis for formulation of the model. As a result, the rates at which population level changes in each compartment become vital. These rates heavily depend on the interaction among the individuals of each compartment. 
\vspace{.2cm} \\

It has been observed that when a disease outbreaks in a population, the healthy individuals tend to change their behavior by adopting protective measures like use of social interventions, pharmaceutical interventions and vaccines etc.. This results to the decrease of rate at which infected individuals grow in the population. So  it  is also necessary for disease models to quantify these interventions.

\newpage

\section{Literature Survey with Reference to COVID-19 Modeling}

\vspace{0.5cm}

Many mathematical models, using ordinary differential equations and delay differential equations, were developed in  to analyze the complex transmission pattern of COVID-19. In \cite{cooper2020sir}, a SIR model is  investigated to study the effectiveness of the modeling approach on the pandemic due to the spread of novel COVID -19 disease. In \cite{ming2020breaking}, a modified SIR epidemic model is stduied to project the actual number of infected cases and the specific burden on isolation wards and intensive care units. In \cite{hernandez2020host}, an in-host modeling study addressed the qualitative characteristics and estimation of standard parameters of coronavirus infection. In \cite{kiselev2021delay,yang2020modeling}, delay differential equations were used  to model the COVID -19 pandemic. Few of the  optimal control studies for COVID-19 involving control interventions such as social interventions, pharmaceutical interventions and vaccines etc. can be found in \cite{aronna2020model, dhaiban2021optimal, kkdjou2020optimal,libotte2020determination, ndondo2021analysis}.
Majorly the mathematical modeling studies in COVID -19 dealt with disease transmission at the population  level. Some limited few studies dealing with spread and control strategies at the age-specific level can be found in \cite{bentout2021age, bubar2021model}. In \cite{hernandez2020host}, an in-host modeling study is discussed on COVID-19. Further the  model parameters were estimated.

\newpage
\section{Objectives}
\vspace{0.5cm}

\begin{itemize}
\item To study the dynamics of COVID-19 with reference to immigration from urban to rural population in Indian scenario.\\

\item To numerically study the stability of equilibria. \\

\item To do the sensitivity analysis for the model parameters. \\

\item To find  parameter regions in where the system is stable using 2-d heat plots. \\ 
\end{itemize}
\newpage

\section{Chapterization}

\vspace{0.5cm}

The chapter-wise division of this work is as follows. In  chapter 2, we formulate the non-linear multi compartmental (SEIR) model and establish the positivity and boundedness of the solution followed by the existence and uniqueness of the solution. In chapter 3, we find the equilibria and derive basic reproduction number of the system. Chapter 4 deals with numerical studies on the  local and global dynamics of the system followed by sensitivity analysis of the model parameters. Later we find the parameter regions in which the system is stable via 2-d heat plots. In  chapter 5  we do the comparative effectiveness studies with reference to different control interventions. Finally in chapter 6 we deal with the discussions and conclusions of the proposed research work.

 \chapter{   SEIR Immigration  Multi Compartment Model }

In this chapter, we initially formulate a non-linear SEIR immigration  multi compartment model and describe the various compartments involved in the model. Later we establish the positivity and boudedness of the solution of the proposed model followed by the existence and uniqueness of the solution.
\newpage
\section{The Mathematical Model And It's Formulation }
We are dividing the entire population into two  groups, urban and rural. We assume that there's an immigration happening from the urban population to rural population as was the case with COVID-19 first wave {\cite{BIF, GLB}}. We also assume that among the immigrants from urban to rural, apart  are quarantined and the remaining directly move to the susceptible of the rural population.
 \\
 
Based on the above assumptions and considerations, we propose the following SEIR immigration  multi compartment model  is given by  the  system of ordinary differential equations: \\

\begin{eqnarray}
   	\frac{dS_{u}}{dt}&=& b_{1} \ -\frac{ \beta S_{u}I_{u}}{N} \ - \mu_{c} S_{u} \ - m S_{u} \label{sec2equ1} \\ 
   	\frac{dE_{u}}{dt}&=& \frac{ \beta S_{u}I_{u}}{N}\ - k E_{u}\ -\mu_{c} E_{u} \ - m E_{u} \label{sec2equ2}\\ 
   	\frac{dI_{u}}{dt} &=&   k {E_{u}} \ -  \gamma I_{u} \ -\mu_{c} I_{u} \label{sec2equ3}\\ 
   	\frac{dR_{u}}{dt} &=& \gamma I_{u}\ - \mu_{c} R_{u} \label{sec2equ4}\\ 
   	\frac{dQ_{r}}{dt} &=& pm(S{u}+E{u})- d_1 Q_{r} \label{sec2equ5}\\
   	\frac{dS_{r}}{dt}&=&  (1-p)m S{u} \ -\frac{ \beta S_{r}I_{r}}{N} \ - \mu_{c} S_{r}  \label{sec2equ6} \\ 
   	\frac{dE_{r}}{dt}&=& (1-p)m E{u} \ + \frac{ \beta S_{r}I_{r}}{N} \ - k E_{r}  \ -  \mu_{c} E_{r}  \label{sec2equ7}\\ 
   	\frac{dI_{u}}{dt} &=&   k {E_{r}} \ -  \gamma I_{r} \ -\mu_{c} I_{r} \label{sec2equ8}\\ 
   	\frac{dR_{r}}{dt} &=& \gamma I_{r}\ - \mu_{c} R_{r} \label{sec2equ9}
   		  \end{eqnarray}

    \begin{table}[ht!]
	\caption{Parameters and their Meanings.} 
     	
     	\centering 
     	\begin{tabular}{|l|l|} 
     		\hline\hline
     		
     		\textbf{Symbols} &  \textbf{Biological Meaning} \\  
     		\hline\hline 
     		$S_u$ & Susceptible urban population  \\
     			\hline\hline 
     		$S_r$ & Susceptible rural population  \\
     		
     		\hline\hline 
     		$E_u$ & Exposed urban population  \\
     			\hline\hline 
     		$E_r$ & Exposed rural population  \\
     			\hline\hline 
     		$I_u$ & Infected urban population  \\
     			\hline\hline 
     		$I_r$ & Infected rural population  \\
     		\hline\hline
     		$R_u$ & Recovered urban population  \\
     			\hline\hline 
     		$R_r$ & Recovered rural population  \\
     		\hline\hline
     		$b_{1}$ & Constant birth rate   \\
     		\hline\hline
     	   	$\beta$ &   transmission rate \\
     		
     		\hline\hline
     		$\mu_{c}$ & natural death rate \\
     	
     		\hline\hline
     		$\gamma$ & recovery rate \\
     	
     		\hline\hline
     		$d_{1}$ & disease induced death rate of population \\
     		 
     		\hline\hline
     $k$ & Incubation  rate 	\\
     \hline\hline	
     		
     $m$ & maturation rate 	\\
     \hline\hline
     	\end{tabular}
     \end{table} \vspace{.25cm}

  \cleardoublepage
  \section{Positivity and Boundedness }

      \underline{\textbf{Positivity of solution}}\textbf{:}   \\ 
     
      We show that when  the initial conditions of the system (2.1.1)-(2.1.9) are positive, then the solution tends to be  positive at  any future time.  A positivity of the solutions are established in similar lines to the method discussed in \cite{kumar2019role, mandale2021dynamics}. Using the  equations (2.1.1)-(2.1.9), we get,

\begin{align*}
\frac{dS_u}{dt} \bigg|_{S_u=0} &=  b_1 \geq 0 ,  &  
\frac{dE_u}{dt} \bigg|_{E_u=0} &= \frac{\beta S_u I_u}{N}\geq 0 ,\\ \\
\frac{dI_u}{dt} \bigg|_{I_u=0} &= k E_u  \geq 0,&
\frac{dR_u}{dt} \bigg|_{R_u=0} &= \gamma I_u \geq 0 , \\ \\
\frac{dQ_r}{dt} \bigg|_{Q_r=0} &= pm(S_u+E_u) \geq 0,\\ \\
\frac{dS_r}{dt} \bigg|_{S_r=0} &= (1-p)mS_u\geq 0 , &
\frac{dE_r}{dt} \bigg|_{E_r=0} &= (1-p)mE_u+ \frac{\beta S_r I_r}{N}\geq 0 ,\\ \\
\frac{dI_r}{dt} \bigg|_{I_r=0} &= kE_r  \geq 0,&
\frac{dR_r}{dt} \bigg|_{R_r=0} &=\gamma I_r \geq 0 & 
\end{align*}

\noindent
\\ 
Thus all the above rates are non-negative on the bounding planes (given by $S_1=0$, $I_1=0$, $R_1=0$, $S_2=0$,  $I_2=0$, and $R_2=0$) of the non-negative region of the real space. So, if a solution begins in the interior of this region, it will remain inside it throughout time $t$. This  happens because the direction of the vector field is always in the inward direction on the bounding planes as indicated by the above inequalities. Hence, we conclude that all the solutions of the the system (2.1.1)-(2.1.9) remain positive for any time $t>0$  provided that the initial conditions are positive. This establishes the positivity of the solutions of the system (2.1.1)-(2.1.9). Next we will show that the solution is bounded.  \vspace{.25cm}
\cleardoublepage

\underline{\textbf{Boundedness of solution}}\textbf{:} \\

Let  $N(t) = S_U(t)+E_U(t))+I_U(t)+R_U(t) +Q_r(t)+ S_r(t)+E_r(t)+I_r(t)+R_r(t) $ \\

Now,  
\begin{equation*}
\begin{split}
\frac{dN}{dt} & = \frac{dS_u}{dt} + \frac{dE_u}{dt} + \frac{dI_u}{dt}+ \frac{d R_u}{dt} +\frac{dQ_r}{dt} +\frac{d S_r}{dt} + \frac{dE_r}{dt}+ \frac{dI_r}{dt}+\frac{d R_r}{dt}  \\[4pt]
& \le b_1-\mu(S_u+E_u+I_u+R_u+Q_r+S_r+E_r+I_r+R_r)  \\
\end{split}
\end{equation*}

The integrating factor here is $e^{\mu t}.$ so, after integration we get, 

$N(t)\le \frac{b_1}{\mu} + ce^{-\mu t}$, as, $c$ is constant. Now  as $t \rightarrow \infty$ we get, 
$$\text{lim sup N(t)} \le \frac{b_1}{\mu}$$

Thus here  we show that the system (2.1.1)-(2.1.9) is positive and bounded.And hence the biologically feasible region is given by the following set, 
\begin{equation*}
\Omega=\bigg \{\bigg(S_u(t), E_u(t), I_u(t), R_u(t), Q_r(t),S_r(t), E_r(t), I_r(t), R_r(t)\bigg) \in \mathbb{R}^{9}_{+}  \\
:(S_u(t) +E_u(t)+ I_u(t)+ R_u(t)+ Q_r(t)+S_r(t)+ E_r(t) +I_r(t)+ R_r(t) \leq \frac{b_1}{\mu}, \ t \geq 0 \bigg\}
\end{equation*}  

  \cleardoublepage
\section{ Existence and Uniqueness of Solutions}

\underline{\textbf{Existence and Uniqueness of Solution}} \\

For the general first order ODE of the form
\begin{equation}
    \dot x=f(t,x) , \hspace{2cm}x(t_0)=x_0
    \end{equation}
    
\noindent
We use the following theorem from \cite{bishal1} in order to establish the existence and uniqueness of solution of the system $(2.1.1)-(2.1.9)$.

\begin{theorem}
Let D denote the domain: 
$$|t-t_0| \leq a, ||x-x_0|| \leq b, x=(x_1, x_2,..., x_n), x_0=(x_{10},..,x_{n0})$$
and suppose that $f(t, x)$ satisfies the Lipschitz condition:
\begin{equation}
||f(t,x_2)-f(t, x_1)|| \leq k ||x_2-x_1||
\end{equation}
and whenever the pairs $(t,x_1)$ and $(t, x_2)$ belong to the domain $D$ , where $k$ is used to represent a positive constant. Then, there exist a constant $\delta > 0$ such that  a unique (exactly one) continuous vector solution $x(t)$ exists for the
system $(2.3.1)$ in the interval $|t-t_0|\leq \delta $. It is important to note that condition $(2.3.2)$ is satisfied by requirement that:
$$\frac{\partial f_i}{\partial x_j}, i,j=1,2,.., n$$
be continuous and bounded in the domain D.
\end{theorem}

\noindent
We use boundedness of the solutions proved above and show that a unique solution exists for system $(2.1.1 )-(2.1.9)$  by showing partial derivative of right hand side of equations $(2.1.1 )-(2.1.9)$ are continuous and bounded with respect to each of the variables $S_u,E_u ,I_u , R_u,Q_r,S_r,E_r  ,I_r$ and $R_r$. \\

Let
 \begin{eqnarray}
 f_1& = &  b_{1} \ - \frac{\beta S_{u}I_{u} }{N} \ - \mu_{c} S_{u} \ - m S_{u} \label{sec2equ11} \\ 
 f_2& = &  \frac{\beta S_{u}I_{u} }{N}\ - k E_{u}  \ -  \mu_{c} E_{u} \ - m E_{u} \label{sec2equ12}\\ 
 f_3& = &     k {E_{u}} \ -  \gamma I_{u} \ -\mu_{c} I_{u} \label{sec2equ13}\\ 
 f_4& = &  	 \gamma I_{u}\ - \mu_{c} R_{u} \label{sec2equ14}\\ 
 f_5& = &  	 pm(S{u}+E{u})- d_1 Q_{r} \label{sec2equ15}\\ 
 f_6& = &  	(1-p)m S{u} \ - \frac{\beta S_{r}I_{r} }{N}  \ - \mu_{c} S_{r}  \label{sec2equ16} \\ 
 f_7& = &  	(1-p)m E{u} \ + \frac{\beta S_{r}I_{r} }{N}\ - k E_{r}  \ -  \mu_{c} E_{r}  \label{sec2equ17}\\ 
 f_8& = &   k {E_{r}} \ -  \gamma I_{r} \ -\mu_{c} I_{r} \label{sec2equ18}\\ 
 f_9& = &  	\gamma I_{r}\ - \mu_{c} R_{r} \label{sec2equ19}
  \end{eqnarray}
   
   \noindent
 From equation $(\ref{sec2equ11})$ we have
 \begin{eqnarray*}
    \frac{\partial f_1}{\partial S_u}&=&-\frac{\beta I_u }{N}, \hspace{.2cm} |\frac{\partial f_1}{\partial S_u}|=|-\frac{\beta I_u }{N}| < \infty\\ \\
    \frac{\partial f_1}{\partial E_u}&=&0, \hspace{.2cm} |\frac{\partial f_1}{\partial E_u}| < \infty\\ \\
 \frac{\partial f_1}{\partial I_u}&=&-\frac{\beta S_u }{N}, \hspace{.2cm}|\frac{\partial f_1}{\partial I_u}| =|-\frac{\beta S_u }{N}| < \infty\\ \\
 \frac{\partial f_1}{\partial R_u}&=&0, \hspace{.2cm}|\frac{\partial f_1}{\partial R_u}|< \infty\\  \\
 \frac{\partial f_1}{\partial Q_r}&=&0, \hspace{.2cm}|\frac{\partial f_1}{\partial Q_r}|< \infty\\ \\
 \frac{\partial f_1}{\partial S_r}&=&0, \hspace{.2cm}|\frac{\partial f_1}{\partial S_r}| < \infty\\ \\
 \frac{\partial f_1}{\partial E_r}&=&0, \hspace{.2cm}|\frac{\partial f_1}{\partial E_r}| < \infty\\ \\
 \frac{\partial f_1}{\partial I_r}&=&0, \hspace{.2cm}|\frac{\partial f_1}{\partial I_r}| < \infty\\ \\
 \frac{\partial f_1}{\partial R_r}&=&0, \hspace{.2cm} |\frac{\partial f_1}{\partial R_r}| < \infty
  \end{eqnarray*}
  \noindent
 From equation $(\ref{sec2equ12})$ we have
 \begin{eqnarray*}
    \frac{\partial f_2}{\partial S_u}&=&\frac{\beta S_u }{N}, \hspace{.2cm} |\frac{\partial f_2}{\partial S_u}|=|\frac{\beta S_u }{N}| < \infty\\ \\
    \frac{\partial f_2}{\partial E_u}&=&-(k+m+\mu_{c}), \hspace{.2cm} |\frac{\partial f_2}{\partial E_u}|= |-(k+m+\mu_{c})| < \infty\\ \\
 \frac{\partial f_2}{\partial I_u}&=&\frac{\beta_1 S_u }{N}, \hspace{.2cm}|\frac{\partial f_2}{\partial I_u}| =|\frac{\beta_1 S_u }{N}| < \infty\\ \\
 \frac{\partial f_2}{\partial R_u}&=&0, \hspace{.2cm}|\frac{\partial f_2}{\partial R_u}|< \infty\\ \\ 
 \frac{\partial f_2}{\partial Q_r}&=&0, \hspace{.2cm}|\frac{\partial f_2}{\partial Q_r}|< \infty\\ \\ 
 \frac{\partial f_2}{\partial S_r}&=&0, \hspace{.2cm}|\frac{\partial f_2}{\partial S_r}| < \infty\\ \\
 \frac{\partial f_2}{\partial E_r}&=&0, \hspace{.2cm}|\frac{\partial f_2}{\partial E_r}| < \infty\\ \\
 \frac{\partial f_2}{\partial I_r}&=&0, \hspace{.2cm}|\frac{\partial f_2}{\partial I_r}|< \infty\\ \\
 \frac{\partial f_2}{\partial R_r}&=&0, \hspace{.2cm} |\frac{\partial f_2}{\partial R_r}| < \infty
  \end{eqnarray*}
  \noindent
 From equation $(\ref{sec2equ13})$ we have
 \begin{eqnarray*}
    \frac{\partial f_3}{\partial S_u}&=&0, \hspace{.2cm} |\frac{\partial f_3}{\partial S_u}| < \infty\\ \\
    \frac{\partial f_3}{\partial E_u}&=& k, \hspace{.2cm} |\frac{\partial f_3}{\partial E_u}|= |k| < \infty\\ \\
 \frac{\partial f_3}{\partial I_u}&=&-(\gamma+\mu_{c}), \hspace{.2cm}|\frac{\partial f_3}{\partial I_u}| =|-(\gamma+\mu_{c})| < \infty\\ \\
 \frac{\partial f_3}{\partial R_u}&=&0, \hspace{.2cm}|\frac{\partial f_3}{\partial R_u}|< \infty\\ \\ 
 \frac{\partial f_3}{\partial Q_r}&=&0, \hspace{.2cm}|\frac{\partial f_3}{\partial Q_r}|< \infty\\ \\ 
 \frac{\partial f_3}{\partial S_r}&=&0, \hspace{.2cm}|\frac{\partial f_3}{\partial S_r}| < \infty\\ \\
 \frac{\partial f_3}{\partial E_r}&=&0, \hspace{.2cm}|\frac{\partial f_3}{\partial E_r}| < \infty\\ \\
 \frac{\partial f_3}{\partial I_r}&=&0, \hspace{.2cm}|\frac{\partial f_3}{\partial I_r}| < \infty\\ \\
 \frac{\partial f_3}{\partial R_r}&=&0, \hspace{.2cm} |\frac{\partial f_2}{\partial R_r}| < \infty
  \end{eqnarray*}
  \noindent
 From equation $(\ref{sec2equ14})$ we have
 \begin{eqnarray*}
    \frac{\partial f_4}{\partial S_u}&=&0, \hspace{.2cm} |\frac{\partial f_4}{\partial S_u}| < \infty\\ \\
    \frac{\partial f_4}{\partial E_u}&=&0, \hspace{.2cm} |\frac{\partial f_4}{\partial E_u}| < \infty\\ \\
 \frac{\partial f_4}{\partial I_u}&=&\gamma, \hspace{.2cm}|\frac{\partial f_4}{\partial I_u}| =|\gamma| < \infty\\ \\
 \frac{\partial f_4}{\partial R_u}&=&-\mu_{c}, \hspace{.2cm}|\frac{\partial f_4}{\partial R_u}|=|-\mu_{c}|< \infty\\ \\ 
 \frac{\partial f_4}{\partial Q_r}&=&0, \hspace{.2cm}|\frac{\partial f_4}{\partial Q_r}|< \infty\\ \\ 
 \frac{\partial f_4}{\partial S_r}&=&0, \hspace{.2cm}|\frac{\partial f_4}{\partial S_r}| < \infty\\ \\
 \frac{\partial f_4}{\partial E_r}&=&0, \hspace{.2cm}|\frac{\partial f_4}{\partial E_r}| < \infty\\ \\
 \frac{\partial f_4}{\partial I_r}&=&0, \hspace{.2cm}|\frac{\partial f_4}{\partial I_r}| < \infty\\ \\
 \frac{\partial f_4}{\partial R_r}&=&0, \hspace{.2cm} |\frac{\partial f_4}{\partial R_r}| < \infty
  \end{eqnarray*}
  \noindent
 From equation $(\ref{sec2equ15})$ we have
 \begin{eqnarray*}
    \frac{\partial f_5}{\partial S_u}&=&pm, \hspace{.2cm} |\frac{\partial f_5}{\partial S_u}|=|pm| < \infty\\ \\
    \frac{\partial f_5}{\partial E_u}&=&pm, \hspace{.2cm} |\frac{\partial f_5}{\partial E_u}|=|pm| < \infty\\ \\
 \frac{\partial f_5}{\partial I_u}&=&-\frac{\beta_1 S_u }{N}|, \hspace{.2cm}|\frac{\partial f_5}{\partial I_u}| =|-\frac{\beta_1 S_u }{N}| < \infty\\ \\
 \frac{\partial f_5}{\partial R_u}&=&0, \hspace{.2cm}|\frac{\partial f_5}{\partial R_u}|< \infty\\ \\ 
 \frac{\partial f_5}{\partial Q_r}&=&-d_{1}, \hspace{.2cm}|\frac{\partial f_5}{\partial Q_r}|=|-d_{1}|< \infty\\ \\ 
 \frac{\partial f_5}{\partial S_r}&=&0, \hspace{.2cm}|\frac{\partial f_5}{\partial S_r}| < \infty\\ \\
 \frac{\partial f_5}{\partial E_r}&=&0, \hspace{.2cm}|\frac{\partial f_5}{\partial E_r}| < \infty\\ \\
 \frac{\partial f_5}{\partial I_r}&=&0, \hspace{.2cm}|\frac{\partial f_5}{\partial I_r} < \infty\\ \\
 \frac{\partial f_5}{\partial R_r}&=&0, \hspace{.2cm} |\frac{\partial f_5}{\partial R_r}| < \infty
  \end{eqnarray*}
  \noindent
 From equation $(\ref{sec2equ16})$ we have
 \begin{eqnarray*}
    \frac{\partial f_6}{\partial S_u}&=&(1-p)m, \hspace{.2cm} |\frac{\partial f_6}{\partial S_u}|=|(1-p)m | < \infty\\ \\
    \frac{\partial f_6}{\partial E_u}&=&0, \hspace{.2cm} |\frac{\partial f_6}{\partial E_u}| < \infty\\ \\
 \frac{\partial f_6}{\partial I_u}&=&-\frac{\beta_1 S_u }{N}, \hspace{.2cm}|\frac{\partial f_6}{\partial I_u}| =|-\frac{\beta_1 S_u }{N}| < \infty\\ \\
 \frac{\partial f_6}{\partial R_u}&=&0, \hspace{.2cm}|\frac{\partial f_6}{\partial R_u}|< \infty\\ \\ 
 \frac{\partial f_6}{\partial Q_r}&=&0, \hspace{.2cm}|\frac{\partial f_6}{\partial Q_r}|< \infty\\ \\ 
 \frac{\partial f_6}{\partial S_r}&=&-\frac{\beta I_r }{N}-\mu_{c}, \hspace{.2cm}|\frac{\partial f_6}{\partial S_r}| =|-\frac{\beta I_r }{N}-\mu_{c}|< \infty\\ \\
 \frac{\partial f_6}{\partial E_r}&=&0, \hspace{.2cm}|\frac{\partial f_6}{\partial E_r}| < \infty\\ \\
 \frac{\partial f_6}{\partial I_r}&=&-\frac{\beta s_r }{N}, \hspace{.2cm}|\frac{\partial f_6}{\partial I_r}|=-\frac{\beta s_r }{N} < \infty\\ \\
 \frac{\partial f_6}{\partial R_r}&=&0, \hspace{.2cm} |\frac{\partial f_6}{\partial R_r}| < \infty
  \end{eqnarray*}
  \noindent
 From equation $(\ref{sec2equ17})$ we have
 \begin{eqnarray*}
    \frac{\partial f_7}{\partial S_u}&=&0, \hspace{.2cm} |\frac{\partial f_7}{\partial S_u}| < \infty\\ \\
    \frac{\partial f_7}{\partial E_u}&=&(1-p)m, \hspace{.2cm} |\frac{\partial f_7}{\partial E_u}|=|(1-p)m| < \infty\\ \\
 \frac{\partial f_7}{\partial I_u}&=&0, \hspace{.2cm}|\frac{\partial f_7}{\partial I_u}|  < \infty\\ \\
 \frac{\partial f_7}{\partial R_u}&=&0, \hspace{.2cm}|\frac{\partial f_7}{\partial R_u}|< \infty\\ \\ 
 \frac{\partial f_7}{\partial Q_r}&=&0, \hspace{.2cm}|\frac{\partial f_7}{\partial Q_r}|< \infty\\ \\ 
 \frac{\partial f_7}{\partial S_r}&=&-\frac{\beta I_r }{N}, \hspace{.2cm}|\frac{\partial f_7}{\partial S_r}|=|-\frac{\beta I_r }{N}| < \infty\\ \\
 \frac{\partial f_7}{\partial E_r}&=&-(k+\mu{c}), \hspace{.2cm}|\frac{\partial f_7}{\partial E_r}|=|-(k+\mu{c})| < \infty\\ \\
 \frac{\partial f_7}{\partial I_r}&=&-\frac{\beta s_r }{N}, \hspace{.2cm}|\frac{\partial f_7}{\partial I_r}|=|-\frac{\beta s_r }{N}| < \infty\\ \\
 \frac{\partial f_7}{\partial R_r}&=&0, \hspace{.2cm} |\frac{\partial f_1}{\partial R_7}| < \infty
  \end{eqnarray*}
  \noindent
 From equation $(\ref{sec2equ18})$ we have
 \begin{eqnarray*}
    \frac{\partial f_8}{\partial S_u}&=&0, \hspace{.2cm} |\frac{\partial f_8}{\partial S_u}|< \infty\\ \\
    \frac{\partial f_8}{\partial E_u}&=&0, \hspace{.2cm} |\frac{\partial f_8}{\partial E_u}| < \infty\\ \\
 \frac{\partial f_8}{\partial I_u}&=&0, \hspace{.2cm}|\frac{\partial f_8}{\partial I_u}| < \infty\\ \\
 \frac{\partial f_8}{\partial R_u}&=&0, \hspace{.2cm}|\frac{\partial f_8}{\partial R_u}|< \infty\\ \\ 
 \frac{\partial f_8}{\partial Q_r}&=&0, \hspace{.2cm}|\frac{\partial f_8}{\partial Q_r}|< \infty\\ \\ 
 \frac{\partial f_8}{\partial S_r}&=&0, \hspace{.2cm}|\frac{\partial f_8}{\partial S_r}| < \infty\\ \\
 \frac{\partial f_8}{\partial E_r}&=&k, \hspace{.2cm}|\frac{\partial f_8}{\partial E_r}| =|k|< \infty\\ \\
 \frac{\partial f_8}{\partial I_r}&=&-(\gamma+\mu_{c}), \hspace{.2cm}|\frac{\partial f_8}{\partial I_r}|= |-(\gamma+\mu_{c})|< \infty\\ \\
 \frac{\partial f_8}{\partial R_r}&=&0, \hspace{.2cm} |\frac{\partial f_8}{\partial R_r}| < \infty
  \end{eqnarray*}
  \noindent
 From equation $(\ref{sec2equ19})$ we have
 \begin{eqnarray*}
    \frac{\partial f_9}{\partial S_u}&=&0, \hspace{.2cm} |\frac{\partial f_9}{\partial S_u}| < \infty\\ \\
    \frac{\partial f_9}{\partial E_u}&=&0, \hspace{.2cm} |\frac{\partial f_9}{\partial E_u}| < \infty\\ \\
 \frac{\partial f_9}{\partial I_u}&=&0, \hspace{.2cm}|\frac{\partial f_9}{\partial I_u}| < \infty\\ \\
 \frac{\partial f_9}{\partial R_u}&=&0, \hspace{.2cm}|\frac{\partial f_9}{\partial R_u}|< \infty\\ \\ 
 \frac{\partial f_9}{\partial Q_r}&=&0, \hspace{.2cm}|\frac{\partial f_9}{\partial Q_r}|< \infty\\ \\ 
 \frac{\partial f_9}{\partial S_r}&=&0, \hspace{.2cm}|\frac{\partial f_9}{\partial S_r}| < \infty\\ \\
 \frac{\partial f_9}{\partial E_r}&=&0, \hspace{.2cm}|\frac{\partial f_9}{\partial E_r}| < \infty\\ \\
 \frac{\partial f_9}{\partial I_r}&=&\gamma, \hspace{.2cm}|\frac{\partial f_9}{\partial I_r}|=|\gamma| < \infty\\ \\
 \frac{\partial f_9}{\partial R_r}&=&-\mu_{c}, \hspace{.2cm} |\frac{\partial f_9}{\partial R_r}|=|-\mu_{c}| < \infty
  \end{eqnarray*}
  \noindent
 
   \noindent
  Hence, we have shown that the partial derivatives of $f=(f_1, f_2,f_3,f_4,f_5,f_6,f_7,f_8,f_9)$ are continuous and bounded. So,from the conclusions of  theorem 2.3.1, there exists a unique solution of system $(2.1.1 )-(2.1.9)$. 
\chapter[Equilibrium Points And Reproduction Number]{\hyperlink{toc}{Equilibrium Points And Reproduction Number}}
\thispagestyle{empty}

In this chapter we briefly discuss about the equilibrium points and the calculate the reproduction number for the system $(2.1.1)-(2.1.9)$. \\

\noindent
We find that the system $(2.1.1)-(2.1.9)$ admits two equilibrium namely the disease free equilibrium and the infected equilibrium. The disease free equilibrium denoted by $E_{0}$ was found to be,\\

     $E_{0} = (S_u^*, S_r^*, Q_r^*,0,0, 0, 0,0,0,0),$  where, \\
          $$ S_{u}^*=\frac{b_1 }{(\mu + m_{c})}$$
     $$ S_{r}^*=\frac{(1-p) b_{1} m}{\mu_{c}(\mu_{c} + m)}$$
     $$ Q_{r}^*=\frac{p b_{1} m}{\ d_1(\mu_{c} + m)},$$
    and the infected equilibrium is denoted by $E_1$ was found to be,\\

     $E_1=(S_{u1}^*,E_{u1}^*,I_{u1}^*,R_{u1}^*,Q_{r1}^*,S_{r1}^*,E_{r1}^*,I_{r1}^*,R_{r1}^*),$ where,
     
     \begin{equation*}
     \begin{split}
     S_{u1}^*&=\frac{(\gamma + \mu_{c})(N(K+m+\mu_{c}))}{k\beta}\\
     E_{u1}^*&=\frac{(\beta S_{u1}^*I_{u1}^*)}{N(\mu_{c}+m+k)}\\
     I_{u1}^*&=\frac{N[b_1K\beta+(\mu_{c}+m)(\gamma+\mu{c}_)(N(\mu_{c}+m+k))]}{\beta}\\
     R_{u1}^*&=\frac{\gamma I_{u1}^*}{\mu_{c}}\\
     Q_{r1}^*&=\frac{pm(S_{u1}^*+E_{u1}^*)}{d_{1}}\\
     S_{r1}^*&=\frac{(1-p)m S_{u1}^*}{\frac{\beta I_{r1}^*}{N} - \mu_{c}}\\
     E_{r1}^*&=\frac{(1-p)m E_{u1}^*+ \frac{\beta S_{r1}^*I_{r1}^*}{N}}{(\mu_{c}+ k)}\\ \\
        &  I_{r1}^* \ \text {is a solution of the cubic equation} \\
        & \frac{(\gamma +\mu_{c}) }{k} N\beta^2 (I_{r1}^*)^3 - N\beta (I_{r1}^*)^2\bigg[\frac{K(1-p)mE_{u1}^*\beta+\mu_{c}(\gamma+\mu_c)N}{k}\bigg]\\
                    &-(1-p)\mu_c m E_{u1}^* N^2\beta I_{r1}^*+ (1-p)mS_{u1}^* = 0  \\ \\
     R_{r1}^*&=\frac{\mu_{c}}{\gamma I_{r1}^*}
       \end{split}
     \end{equation*} \\

    \noindent 
     From the Descrate's rule of signs, we see that the cubic equation in $I_{r1}^*$ admits a maximum of two positive roots as there are only two sign changes the first sign  change being from first and second terms and the second one being from third to fourth terms. \\
     
        \noindent 
     So the system (2.1.1) - (2.1.9) at most can admit two infected equilibria.
    
      \cleardoublepage

     \underline{Calculation of ${\mathcal R}_0$}
    The basic reproduction number is one of the most important quantities in disease modelling. It is defined as the average number of secondary cases generated for every primary case generated.  \\
    
        \noindent 
    For our proposed model $(2.1.1)-(2.1.9),$ we  calculate the  reproduction number using the next generation matrix technique \cite{diekmann2010construction}.  \\
    
        \noindent 
    As part of this method,  we divide the system is  into four infected and five non-infected states $(2.1.1)-(2.1.9)$. We later obtain the Jacobian matrix of the system's infected states $(2.1.1)-(2.1.9)$ at disease-free equilibrium $E_0$ by calculating the Jacobian matrix of the system's infected states $(2.1.1)-(2.1.9)$ at disease-free equilibrium $E_0$ given by
    
\begin{equation*}
J(E_0)=    
\begin{bmatrix}
\ -k -m -\mu_{c} & 0 & p_1\beta & 0  \\[9pt]
\ (1-p) m  & -k-\mu_{c} & 0 & p_2\beta\\[9pt]
\ k & 0 & -\gamma -\mu_{c} & 0 \\[9pt]
\ 0 & k & 0 & -\gamma -\mu_{c}\\[9pt]
\end{bmatrix}
\end{equation*}      
or, \\

$J(E_0) = T + \sum,$ where, \\

\begin{equation*}
T =    
\begin{bmatrix}
0&0&p_1\beta& 0\\[9pt]
0&0&0&p_2\beta  \\[9pt]
0&0&0&0\\[9pt]
0&0&0&0\\[9pt]
\end{bmatrix}
\end{equation*}

\begin{equation*}
\sum =    
\begin{bmatrix}
-k -m -\mu_{c} &0&0&0\\[9pt]
(1-p)m&-k -m -\mu_{c} &0&0\\[9pt]
k  &0&-\gamma -\mu_{c}&0\\[9pt]
0&k  &0&-\gamma -\mu_{c}\\[9pt]
\end{bmatrix}
\end{equation*}

Calculating the inverse of  $\sum,$ we get,
\begin{equation*}
{\sum}^ {-1} =    
\begin{bmatrix}
\frac{1}{(-k-m- \mu_{c})} & 0 &0&0 \\[9pt]
\frac{(1-p)}{(-k-m- \mu_{c})(-k-\mu_{c}} & \frac{1}{(-k- \mu_{c})} &0&0 \\[9pt]
\frac{(k)}{(k+m+ \mu_{c})(-\gamma-\mu_{c}}&0&\frac{(1)}{(-\gamma-\mu_{c})}&0\\[9pt]
\frac{k(1-p)m}{(-k-m- \mu_{c})(-k-\mu_{c})(\gamma+\mu_{c})}&\frac{k}{(k+ \mu_{c})(-\gamma-\mu_{c}}&0&\frac{(1)}{(-\gamma-\mu_{c})}\\[9pt]
\end{bmatrix}
\end{equation*}
Now
\begin{equation*}
-T*{\sum}^ {-1} =    
\begin{bmatrix}
\frac{\beta p_1 k}{(k+m +\mu_{c})(\gamma+\mu_{c})}&0 & \frac{\beta p_1}{(\gamma+\mu_{c})}&0 \\[9pt]
\frac{\beta p_2(1-p)m k}{(k+m +\mu_{c})(\gamma+\mu_{c})(k+\mu_{c})} & \frac{\beta p_2k}{(\gamma+\mu_{c})(k+\mu_{c})}&0 &\frac{\beta p_2}{(\gamma+\mu_{c})}\\[9pt]
0&0&0&0\\[9pt]
0&0&0&0\\[9pt]
\end{bmatrix}
\end{equation*}
Now
\begin{equation*}
k= E*-T*{\sum}^ {-1}*E^{T}  \ =   \ 
\begin{bmatrix}
\frac{\beta p_1k}{(\gamma+\mu_{c})(k+m+\mu_{c})}&0\\[9pt]
0&\frac{\beta p_2k}{(\gamma+\mu_{c})(k+\mu_{c})}\\[9pt]
\end{bmatrix}
\end{equation*}
 

As per the next generation matrix methos we have,  \\ \\

  {\bf{ $R_0 = { {\frac{\beta p_1k}{(\gamma+\mu_{c})(k+m+\mu_{c})}}},$}} the most dominant eigen value of the matrix $k$.\\ \\ 
\chapter[Numerical Studies]{\hyperlink{toc}{Numerical Studies}}

In this chapter we do the various numerical simulations dealing with the local and global stabilities of  disease free equilibrium and local stability of infected equilibrium. We also perform the  sensitivity analysis of the model parameters for identifying the  the sensitive parameters and the corresponding range. 2-d heat plots are also done for identifying the  parameter regions in which the system is stable.

\newpage

\section{Parameters Values}
    A parameter is a variable that affects the output or behaviour of a mathematical entity yet is considered constant. Parameters and variables are closely connected, and the distinction is sometimes only a question of perspective. Variables are thought to change, whereas parameters are thought to stay the same or change slowly. In certain cases, it's possible to envisage doing several tests with the variables changing from one to the next, but the parameters remaining constant throughout and changing only between trials.

    Table 2 lists the parameter values as well as the source from which they were obtained. With these parameters, the asymptotic stability of will be quantitatively demonstrated in $E_0$ and $E_1$ depending on the values of the fundamental reproduction number in a manner similar to \cite{kouokam2013disease}. 

    \vspace{1cm}
    \begin{table}[ht!]
	\caption{values of parametres and their Source.} 
     	
     	\centering 
     	\begin{tabular}{|l|l|l|} 
     		\hline\hline
     		
     		\textbf{Parameters} &  \textbf{Values}&  \textbf{Source} \\  
     		\hline\hline 
     	
     		$b_{1}$ & $\mu N(0)$ & \cite{samui2020mathematical} \\
     		\hline\hline
     		$\gamma$ & 0.0714 & \cite{kouokam2013disease}  \\
     		\hline\hline
     		$\beta$ &  0.00028 & \cite{srivastav2021modeling}\\
     			\hline\hline
     		$\ d_{1}$ &  0.013 & \cite{srivastav2021modeling}\\	
     		\hline\hline
     		$\mu_{c}$ & 0.0062 & \cite{samui2020mathematical}\\
     	
     		\hline\hline
     	$ k$ & 0.1961 & \cite{kumar2019role}\\
     	\hline\hline
     $m$ & 0.000182 	& 	\cite{kouokam2013disease}\\
     \hline\hline
     	\end{tabular}
     \end{table} \vspace{.25cm}

  \newpage
\section{Numerical Simulations for Disease Free Equilibrium}
\subsection{Local Stability}
\vspace{.50cm}

     We prove numerically that the disease-free equilibrium $E_0$ is locally asymptotically stable whenever ${\mathcal R}_0 <1$. We adjust some of the table 4.2 parameter values to make the value of ${\mathcal R}_0$ smaller than one. When the values of $\beta$, and $\mu$ were set to 0.00028 and 0.62, respectively, the value of ${\mathcal R}_0$ was estimated to be $0.15$ and $E_0 = (487.05, 325.24,45.29,0,0,0,0,0,0) $. The system of equations $(2.1.1)-(2.1.9)$ was numerically solved in MATLAB sofware using parameter values from table 4.2. Figure 4.1 depicts the system solutions $(2.1.1)-(2.1.9)$ with the starting variables $(S_u^*,E_u^*,I_u^*, R_u^*,Q _r^*, S_r^*,E _r^*, I_r^*,R _r^*)$=(100,85,50,20,10,100,85,50,20).\newline Figure 4.1 shows that the solution finally approaches the infection-free condition, $E_0$. As a result of our numerical study, we find that the infection-free equilibrium $E_0$ of the system$(2.1.1)-(2.1.9)$ is locally asymptotically stable when ${\mathcal R}_0$ is smaller than unity. Table 4.2 deals with the  all the parameter values.
     
     \vspace{1cm}
    \begin{table}[ht!]
	\caption{Values of parameters for $ {\mathcal R}_0 < 1 $.} 
     	
     	\centering 
     	\begin{tabular}{|l|l|} 
     		\hline\hline
     		
     		\textbf{Parameters} &  \textbf{Values} \\  
     		\hline\hline 
     	
     		$b_{1}$ & $350$   \\
     		\hline\hline
     		$\gamma$ & 0.0714    \\
     		\hline\hline
     		$\beta$ &  0.00028  \\
     			\hline\hline
     		$\ d_{1}$ &  0.013  \\	
     		\hline\hline
     		$\mu_{c}$ & 0.0062  \\
     	
     		\hline\hline
     	$ k$ & 0.1961  \\
     	\hline\hline
     $m$ & 0.000182 	 	\\
     \hline\hline
     	\end{tabular}
     \end{table} \vspace{.25cm}

     \begin{center}

	\begin{figure}[hbt!]
		\includegraphics[height = 9cm, width = 16 cm]{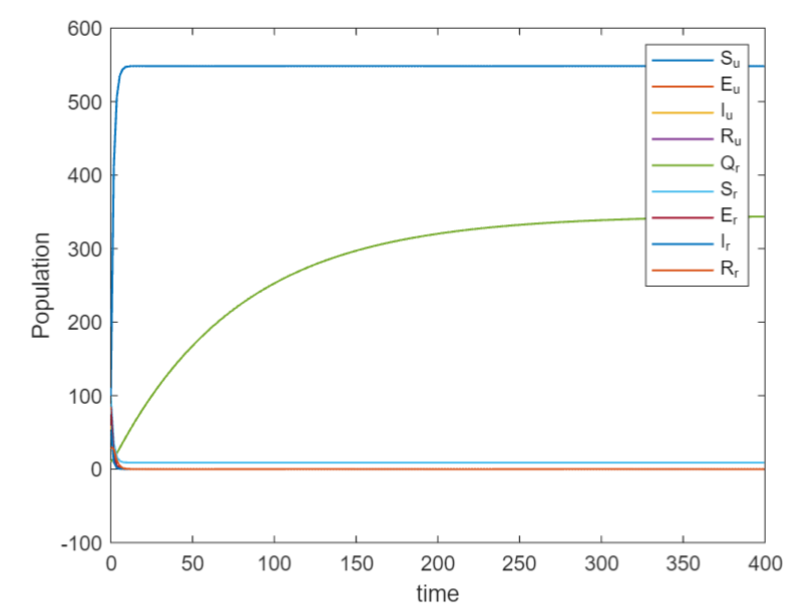} 
			\caption{Figure shows that local asymptotic stability of $E_0$ as $R_0 <1$.}
	
	\end{figure} 

\end{center}

\cleardoublepage

\subsection{Global stability}

For numerically establishing the global stability of disease free equilibrium, we have studies the state trajectories of both urban and rural populations  at different points and observed their convergence. The different initial points for urban populations include   $\{(20,40,10),(107,36,21), (238,46,32), (175,14,30),$ \\ $(175,34,50)\}$         and the corresponding trajectories are depicted in the figure 4.2 and the different initial points for rural populations include   $ \{(50,40,4), $\\ $(63.376,49.1008,30.2826),  (69.376,19.1008,49.2826), (73.654,20.527,9.019),$\\$(80.654,39.527,9.019)\} $      and the corresponding trajectories are depicted in the figure 4.3.

\vspace{8cm}

 \begin{center}
	\begin{figure}[hbt!]
		\includegraphics[height = 9cm, width = 17.5cm]{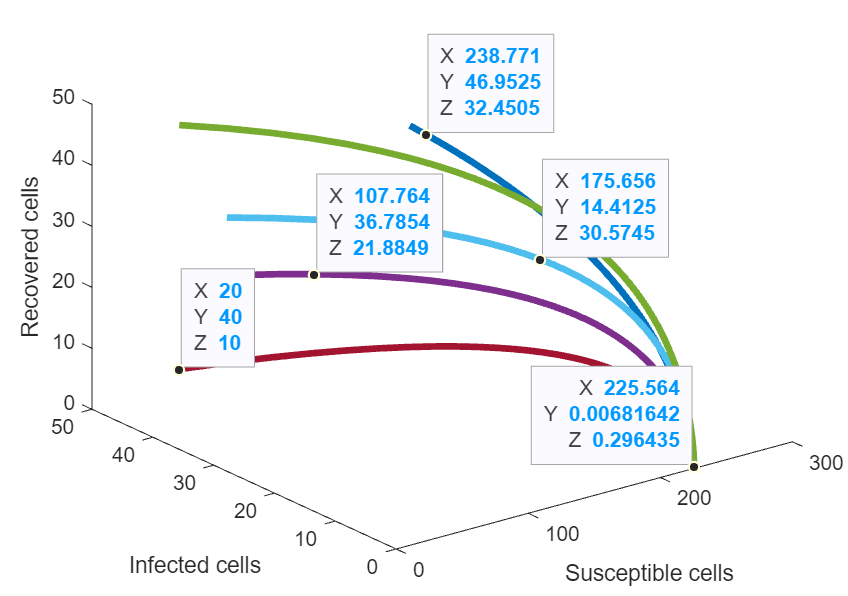} 
			\caption{Figure shows that global asymptotic stability of $E_0$ as $R_0 <1$. in urban scenario.}
	
	\end{figure} 
\end{center}
\newpage 
\begin{center}
	\begin{figure}[hbt!]
		\includegraphics[height = 9cm, width = 17.5cm]{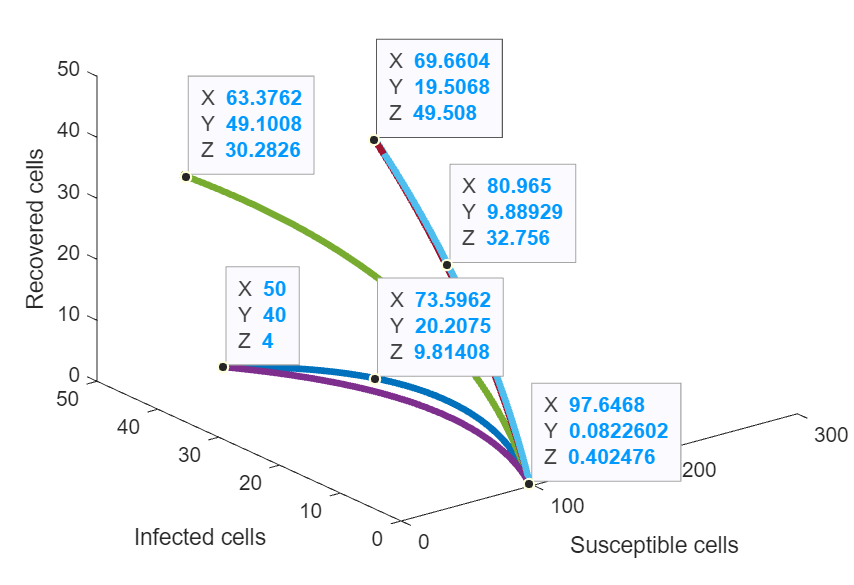} 
			\caption{Figure show that global asymptotic stability of $E_1$ whenever $R_0 < 1$ in rural scenario  .}
	
	\end{figure} 
\end{center}

\cleardoublepage

\section{Stability of infected equilibrium}
\subsection{Local stability}
Though from chapter 3, we see that the system (2.1.1) - (2.1.9), can admit two infected equilibria, from the numerical studies we found the system admits one infected equilibrium that is locally asymptotically stable. The numerical depiction of the same is done below. \\
\vspace{1cm}
    \begin{table}[ht!]
	\caption{Values of parametres for $R_0 > 1.$} 
     	
     	\centering 
     	\begin{tabular}{|l|l|} 
     		\hline\hline
     		
     		\textbf{Parameters} &  \textbf{Values} \\  
     		\hline\hline 
     	
     		$b_{1}$ & $350$   \\
     		\hline\hline
     		$\gamma$ & 0.0714    \\
     		\hline\hline
     		$\beta$ &  0.00028  \\
     			\hline\hline
     		$\ d_{1}$ &  0.013  \\	
     		\hline\hline
     		$\mu_{c}$ & 0.0062  \\
     	
     		\hline\hline
     	$ k$ & 0.1961  \\
     	\hline\hline
     $m$ & 0.000182 	 	\\
     \hline\hline
     	\end{tabular}
     \end{table} \vspace{.25cm}
     \\
    \newpage
    \noindent	

 For the parameter values in table 4.3, the value of ${\mathcal R}_0$ was estimated to be $1.514 $   and the infected equilibrium to be  \\
 $E_1  = (4364.05, 24.45, 35.04, 22.35, 2705.25, 704.48, 45.22, 37.53, 28.21) $.  \\
 
 \noindent
 Figure 4.4 shows that $E_1$ is locally asymptotically stable for $R_0 > 1.$  \\
 
  \noindent
 The initial values for this simulation were chosen to be \\ 
 
 $(S_{u1}^*,E_{u1}^*,I_{u1}^*,R_{u1}^*,Q_{u1}^*,S_{r1}^*,E_{r1}^*,I_{r1}^*,R_{r1}^*) = (100,85,50,20,10,100,85,50,20).$ \\
 
	\begin{figure}[hbt!]
		\includegraphics[height = 9cm, width = 15cm]{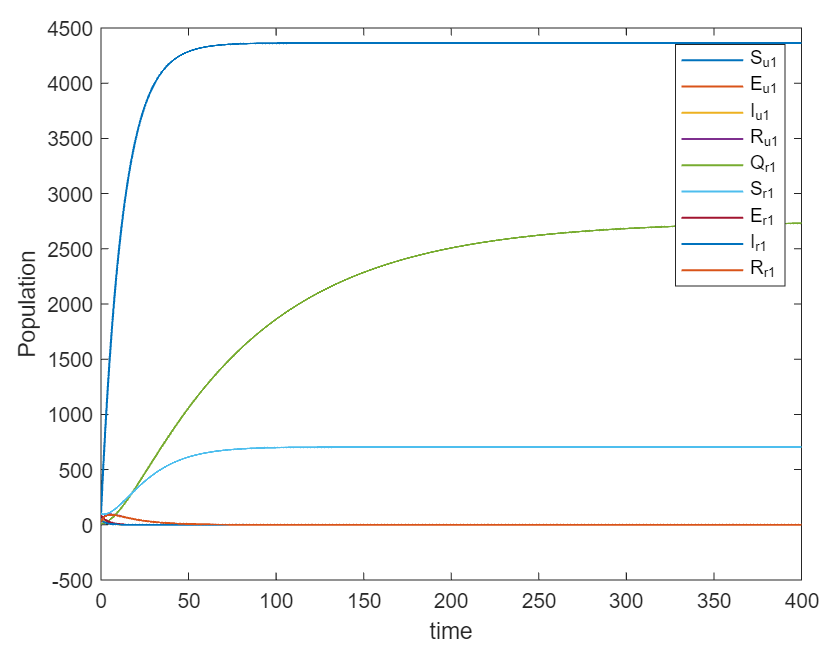} 
			\caption{Figure show that local asymptotic stability of $E_1$ whenever $R_0 > 1$.}
	
	\end{figure} 
	\newpage

\noindent	

 For the parameter values in table 4.3, the value of ${\mathcal R}_0$ was estimated to be $1.514 $   and the infected equilibrium to be  \\
 $E_1 = (4364.05, 24.45, 35.04, 22.35, 2705.25, 704.48, 45.22, 37.53, 28.21)$.  \\
 
 \noindent
 Figure 4.5 shows that $E_1$ is locally asymptotically stable for $R_0 > 1.$  \\
 
  \noindent
 The initial values for this simulation were chosen to be \\
 $(S_{u1}^*,E_{u1}^*,I_{u1}^*,R_{u1}^*,Q_{u1}^*,S_{r1}^*,E_{r1}^*,I_{r1}^*,R_{r1}^*)= (150,85,100,20,10,100,85,80,20).$ \\ 

	\begin{figure}[hbt!]
		\includegraphics[height = 9cm, width = 15cm]{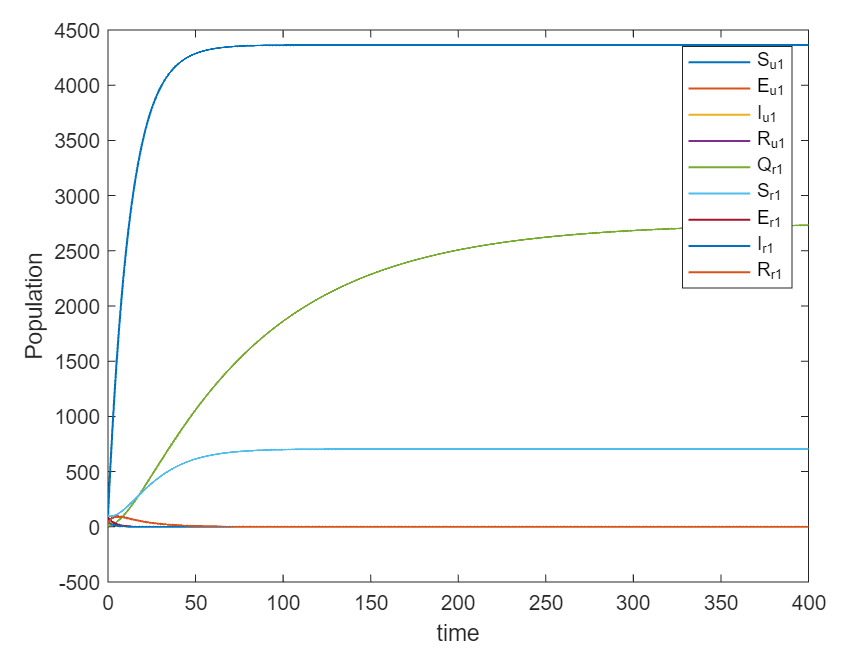} 
			\caption{Figure show that local asymptotic stability of $E_1$ whenever $R_0 > 1$.}
	
	\end{figure} 
 \newpage
\noindent	

 For the parameter values in table 4.3, the value of ${\mathcal R}_0$ was estimated to be $1.514 $   and the infected equilibrium to be  \\
 $E_1 =(4364.05, 24.45, 35.04, 22.35, 2705.25, 704.48, 45.22, 37.53, 28.21) $.  \\
 
 \noindent
 Figure 4.6 shows that $E_1$ is locally asymptotically stable for $R_0 > 1.$  \\
 
  \noindent
 The initial values for this simulation were chosen to be \\  
 
 $(S_{u1}^*,E_{u1}^*,I_{u1}^*,R_{u1}^*,Q_{u1}^*,S_{r1}^*,E_{r1}^*,I_{r1}^*,R_{r1}^*) = (90,85,60,20,10,100,85,70,20).$ \\
 
	\begin{figure}[hbt!]
		\includegraphics[height = 9cm, width = 15cm]{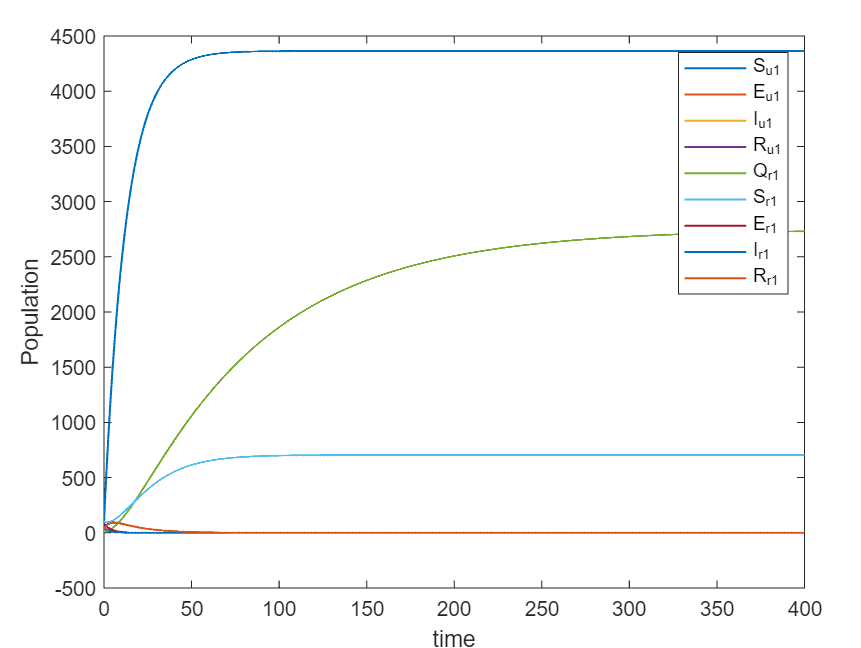} 
			\caption{Figure show that local asymptotic stability of $E_1$ whenever $R_0 > 1$.}
	
	\end{figure} 
 \newpage
\noindent	

 For the parameter values in table 4.3, the value of ${\mathcal R}_0$ was estimated to be $1.514 $   and the infected equilibrium to be  \\
 $E_1 =(4364.05, 24.45, 35.04, 22.35, 2705.25, 704.48, 45.22, 37.53, 28.21) $.  \\
 
 \noindent
 Figure 4.7 shows that $E_1$ is locally asymptotically stable for $R_0 > 1.$  \\
 
  \noindent
 The initial values for this simulation were chosen to be \\
 
 $(S_{u1}^*,E_{u1}^*,I_{u1}^*,R_{u1}^*,Q_{u1}^*,S_{r1}^*,E_{r1}^*,I_{r1}^*,R_{r1}^*) = (200,95,50,40,10,100,85,80,20).$ \\
 \newline
	\begin{figure}[hbt!]
		\includegraphics[height = 9cm, width = 15cm]{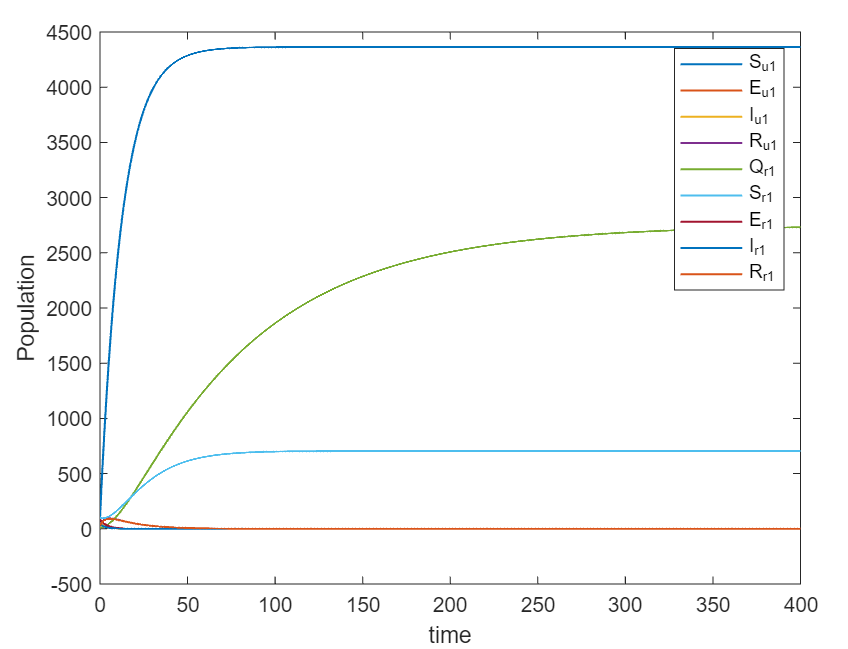} 
			\caption{Figure show that local asymptotic stability of $E_1$ whenever $R_0 > 1$.}
	
	\end{figure} 

\cleardoublepage

\newpage
\section{Sensitivity Analysis}
One of the  greatest concerns in a pandemic is the ability of an infection to infiltrate in a population. Sensitivity analysis is used to investigate the elements that contribute to the spread and persistence of this disease in the community. We are interested in the characteristics that cause a higher divergence in the value of the  reproduction number.\newline

In this setting sensitivity analysis quite handy as it is used to investigate the effect of input factors (parameters, in-output boundary conditions, and so on) on output variables. 

\noindent
The infection dies out in the population when $R _0 <1$, as shown in the preceding sections. As a result, it's critical to keep the model parameters under control such that $R _0$ is smaller than one. As a result, it's critical to figure out what intervals the model parameters are sensitive too. In this part, we do a sensitivity analysis of the model parameters, similar to the sensitivity analysis performed. We depict the infected population, the mean infected population, and the mean squared error as a function of time as each parameter is changed. These charts may be used to see if a parameter is sensitive within a specific interval. The various intervals used are listed in the table \ref{t} below. Table 4.3 contains the fixed parameter values.

\begin{table}[hbt!]
		\caption{Interval Ranges for Sensitivity Analysis}
		\centering
		\label{t}
		{
			\begin{tabular}{|l|l|l|}
				\hline
				\textbf{Parameter} & \textbf{Interval} & \textbf{Step Size}  \\
				\hline

				 $b_1$ & 345 to 355 & 1 
					\\ \cline{2-2}
				 & 355  to 365 &
				 	\\ \cline{2-2}
				 
				 \hline
				 	$ m $ & 0 to 0.00182 & 0.0001 
					\\ \cline{2-3}
				 & 0.00182 to 0.1  & 0.0001 \\
				 \hline

				$\beta$ & 0 to 0.00028 & 0.0001 
				\\ \cline{2-2}
				& 0.00028 to 0.1 &  \\ 
				\hline 
					
				$k$ & 0 to 0.05 & 0.01  
				\\ \cline{2-2}
				& 0.2to 2 & \\
				\hline 
				$d_1$ & 0 to 0.013  & 0.001  
				\\ \cline{2-3}
				& 0.013 to 0.5 & 0.001 \\
				\hline 
			
				$\mu_{c}$ & 0.1 to .5 & .001  
				\\ \cline{2-2}
				& .5 to 1 & \\
				\hline
			$\gamma$ & 0 to 0.0714 & 0.001  
				\\ \cline{2-2}
				& 0.0714 to 1 & \\
				\hline 
			
			\end{tabular}
		}
	\end{table}
	
	\newpage
	
	\underline{{Parameter $\boldsymbol{\mu_{c}}$}}
	
Figure 4.8 depicts the results of the sensitivity of $\mu_c$, which was determined in two intervals as specified in table \ref t. The sensitivity is obtained by plotting the infected population for each value of the parameter $\mu_c$ in the intervals, the mean infected population, and the mean square error. Figure 4.8shows that the overall infected population remains constant for all values of $\mu_c$ changed in both intervals given in table \ref t. The mean infection decreases rapidly and is resolved in a few days. It is evident that the mean square error of the overall infected population increases initially. This variation around the mean, however, lasts only a short time before the mean square error converges to zero. Because there is just one expected variation and the standard deviation falls to minimal levels, we may conclude that $\mu_c$ is insensitive in both intervals I and II.
		
\begin{figure}[hbt!]
\begin{center}
\includegraphics[width=2.2in, height=1.8in, angle=0]{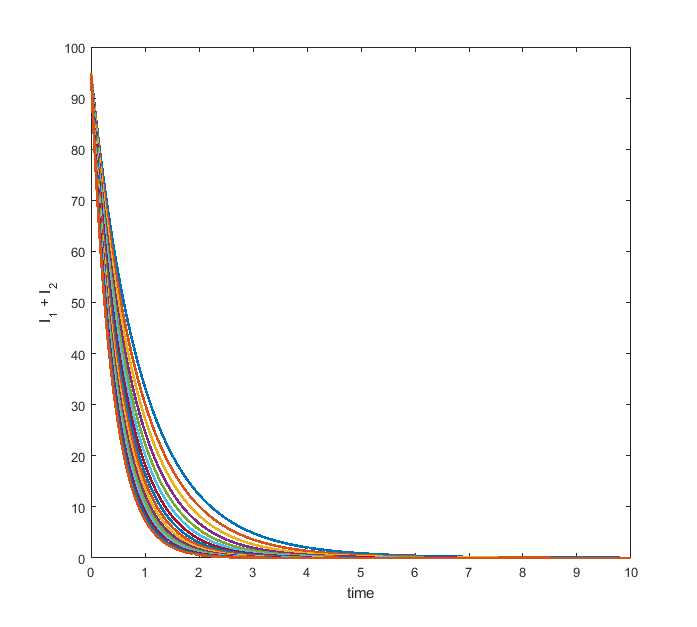}
\hspace{-.4cm}
\includegraphics[width=2.2in, height=1.8in, angle=0]{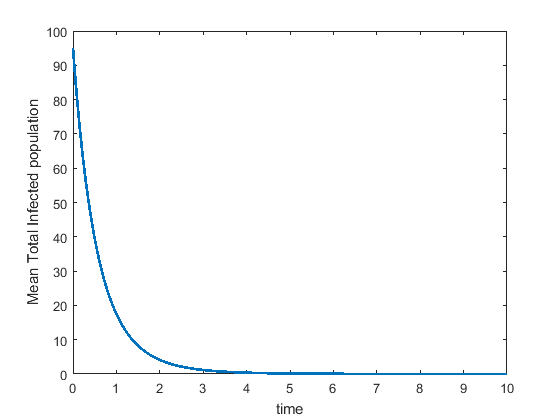}
\hspace{-.395cm}
\includegraphics[width=2.2in, height=1.8in, angle=0]{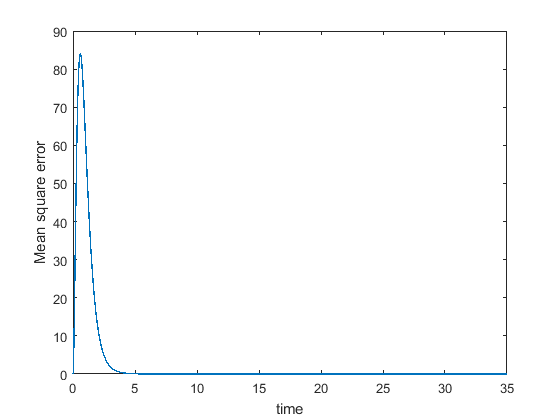}
\caption*{(a) Interval I : 0 to 0.5}
\end{center}
\end{figure}

\vspace{-3mm}

\begin{figure}[hbt!]
\begin{center}
\includegraphics[width=2.2in, height=1.8in, angle=0]{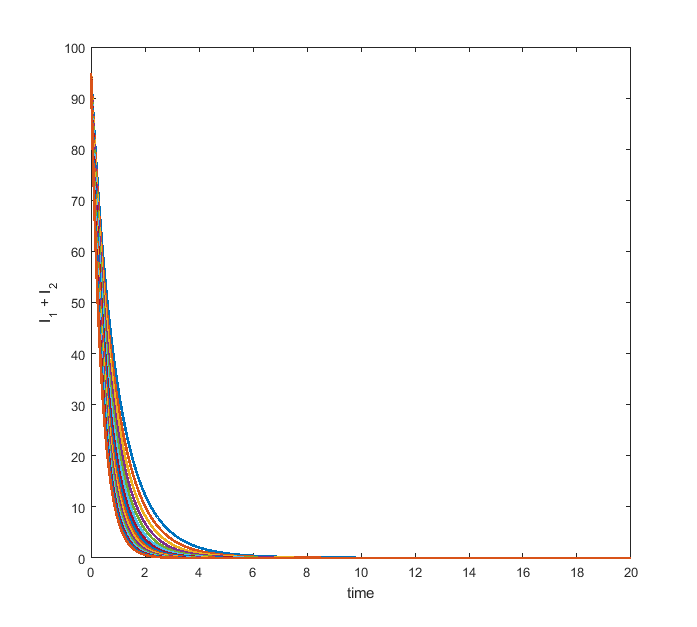}
\hspace{-.4cm}
\includegraphics[width=2.2in, height=1.8in, angle=0]{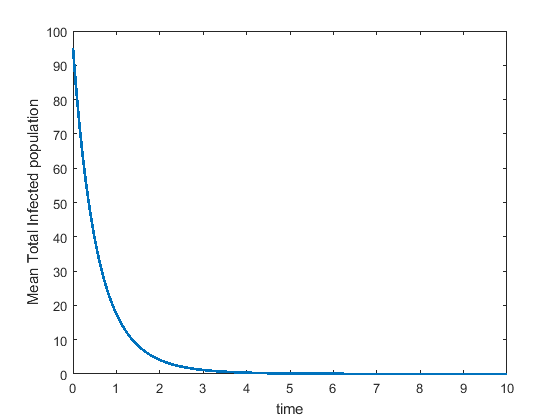}
\hspace{-.395cm}
\includegraphics[width=2.2in, height=1.8in, angle=0]{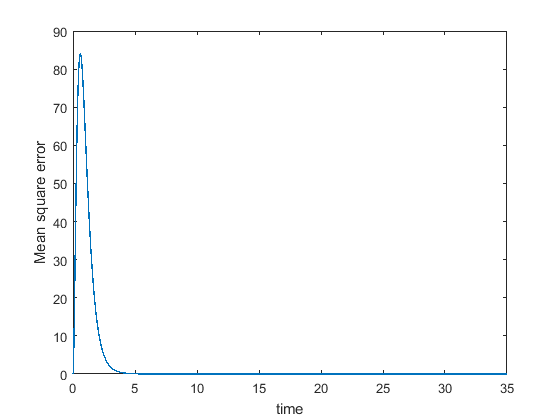}
\caption*{(b) Interval II : 0.5 to 1}
\vspace{5mm}
\caption{Figure showing that the sensitivity analysis of $\mu_{c}$ varied in 2 intervals in table \ref{t}. The plots shows that  the infected population for each  value of the parameter $\mu_{c}$  per every  interval and  with  the mean infected population and the mean square error in the same interval. }
\end{center}
\end{figure}

\newpage
\subsection{{Parameter k}}

The results related to sensitivity of $k$, varied in two intervals as mentioned in table \ref{t}, are given in figure \ref{k}. The plots of infected population for each varied value of the parameter $k$ per interval, the mean infected population and the mean square error are used to determine the sensitivity. We conclude from these plots that the parameter $k$ is sensitive in interval I and  II. 

\newpage
		\begin{figure}[hbt!]
			\begin{center}
				\includegraphics[width=2.2in, height=1.6in, angle=0]{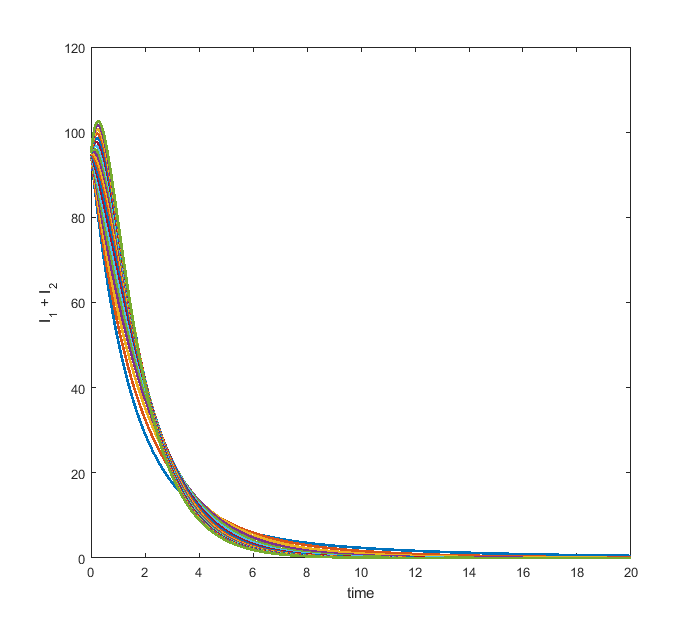}
				\hspace{-.4cm}
				\includegraphics[width=2.2in, height=1.6in, angle=0]{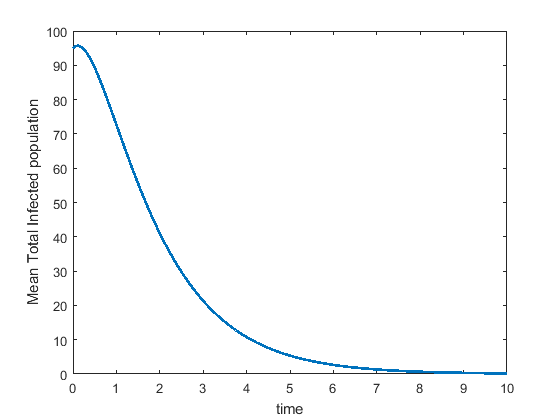}
				\hspace{-.395cm}
					\includegraphics[width=2.2in, height=1.6in, angle=0]{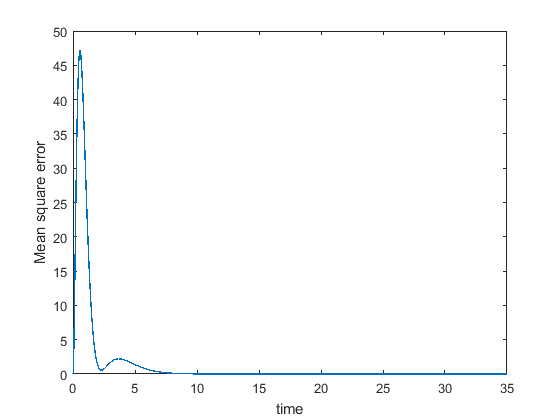}
			\caption*{(a) Interval I :  0 to 1.33}
				
			\end{center}
		\end{figure}
	\begin{figure}[hbt!]
			\begin{center}
				\includegraphics[width=2.2in, height=1.6in, angle=0]{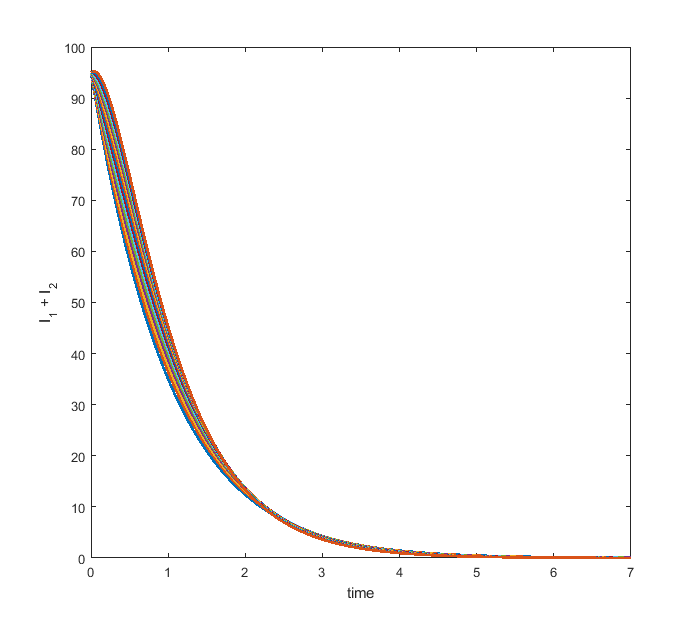}
				\hspace{-.4cm}
				\includegraphics[width=2.2in, height=1.6in, angle=0]{k_2n.png}
				\hspace{-.395cm}
					\includegraphics[width=2.2in, height=1.6in, angle=0]{k_1n.png}
			\caption*{(b) Interval I :  1.33 to 2}

			\caption{Figure showing that the sensitivity analysis of $k$ varied in 2 intervals in table \ref{t}. The plots shows that  the infected population for each  value of the parameter $k$  per every  interval and  with  the mean infected population and the mean square error in the same interval.   }
				\label{k}
			\end{center}
\end{figure}

\newpage
\subsection{{Parameter }}
\underline{$\gamma$}
The results related to sensitivity of $\gamma$, varied in two intervals as mentioned in table \ref {t}, are given in figure \ref {gamma}. The plots of infected population for each varied value of the parameter $\gamma$ per interval, the mean infected population and the mean square error are used to determine the sensitivity. We conclude from these plots that the parameter $\gamma$ is sensitive in interval I and  II. In similar lines, the sensitivity analysis is done for other parameters. The results are summarized below. 
\newpage
		\begin{figure}[hbt!]
			\begin{center}
				\includegraphics[width=2.2in, height=1.6in, angle=0]{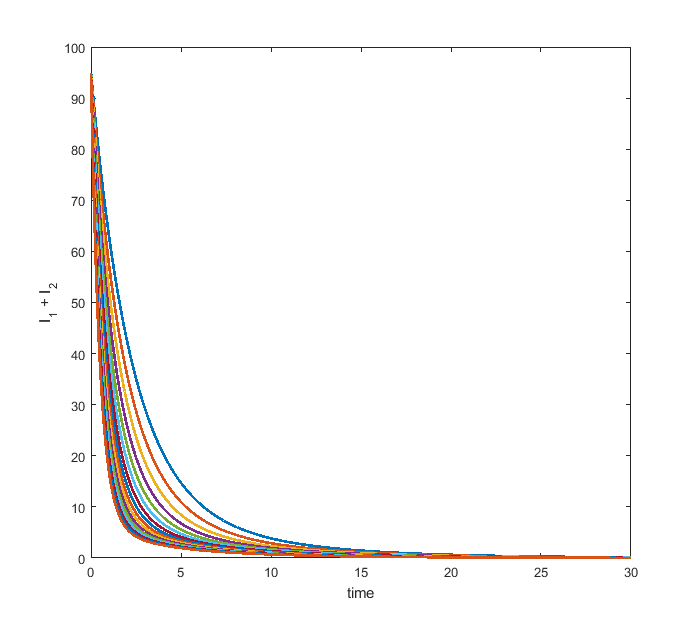}
				\hspace{-.4cm}
				\includegraphics[width=2.2in, height=1.6in, angle=0]{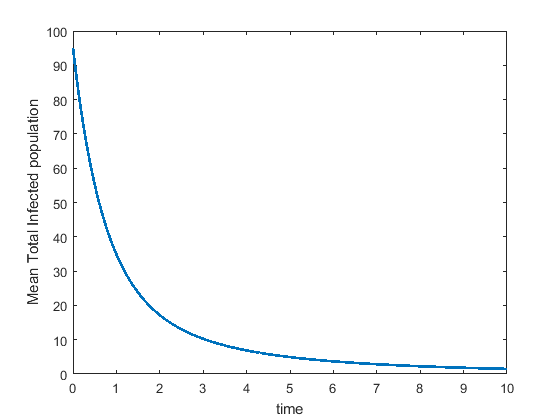}
				\hspace{-.395cm}
					\includegraphics[width=2.2in, height=1.6in, angle=0]{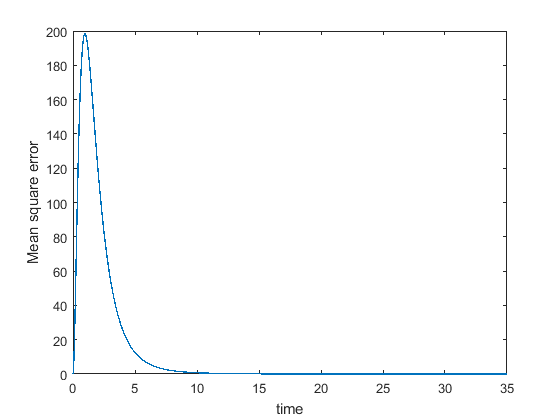}
			\caption*{(a) Interval I :  0 to 1.33}
				
			\end{center}
		\end{figure}
	\begin{figure}[hbt!]
			\begin{center}
				\includegraphics[width=2.2in, height=1.6in, angle=0]{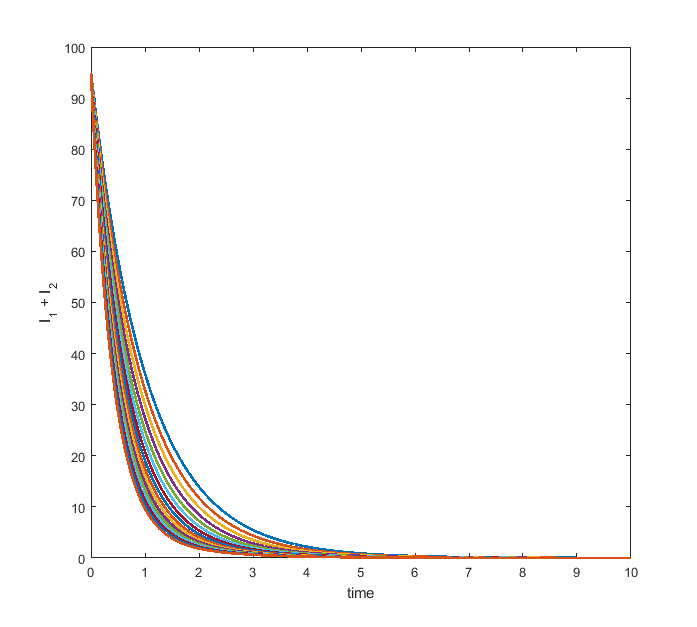}
				\hspace{-.4cm}
				\includegraphics[width=2.2in, height=1.6in, angle=0]{gama.2n.png}
				\hspace{-.395cm}
					\includegraphics[width=2.2in, height=1.6in, angle=0]{gama.3n.png}
			\caption*{(b) Interval I :  1.33 to 2}

			\vspace{7mm}
			\caption{Figure showing that the sensitivity analysis of $\gamma$ varied in 2 intervals in table \ref{t}. The plots shows that  the infected population for each  value of the parameter $\gamma$  per every  interval and  with  the mean infected population and the mean square error in the same interval.   }
				\label{gamma}
			\end{center}
\end{figure}

\cleardoublepage

\underline{{Parameter $\boldsymbol \beta$}}
The results related to sensitivity of $\beta$, varied in two intervals as mentioned in table \ref {t}, are given in figure 4.11. The plots of infected population for each varied value of the parameter $\beta$ per interval, the mean infected population and the mean square error are used to determine the sensitivity. We conclude from these plots that the parameter $\beta$ is insensitive in interval I and  II. In similar lines, the sensitivity analysis is done for other parameters. The results are summarized below. 
\begin{figure}[hbt!]
\begin{center}
\includegraphics[width=2.2in, height=1.8in, angle=0]{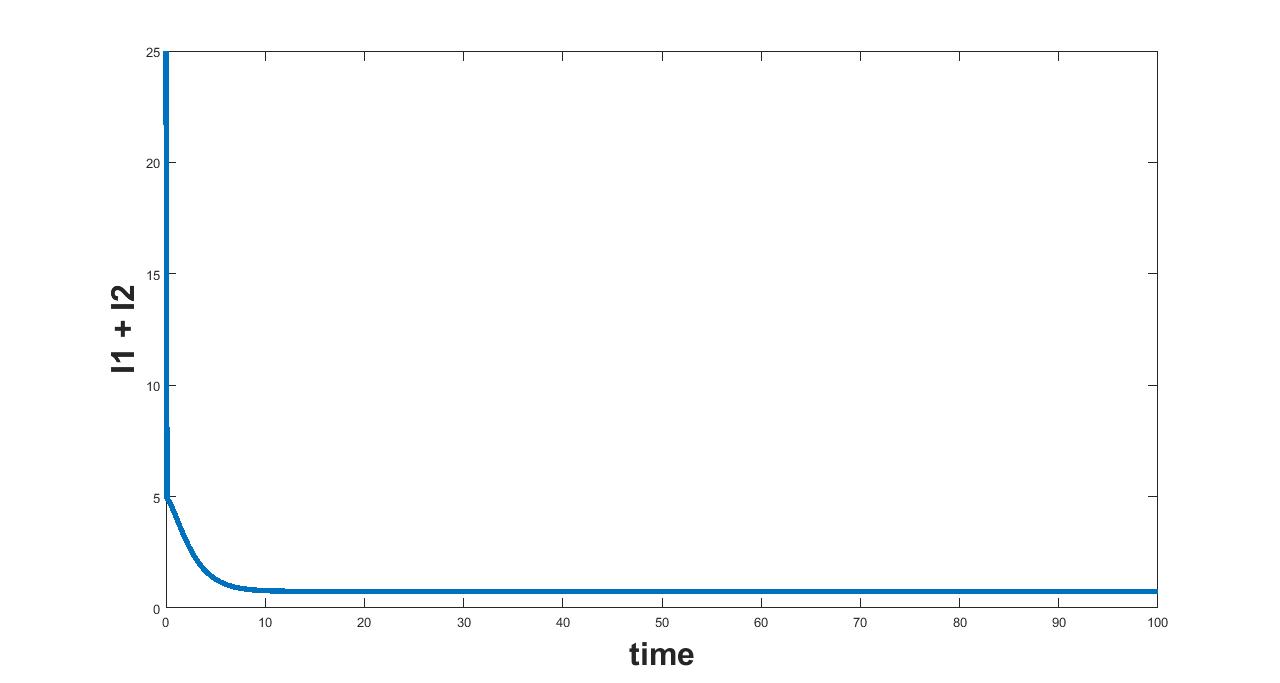}
\hspace{-.4cm}
\includegraphics[width=2.2in, height=1.8in, angle=0]{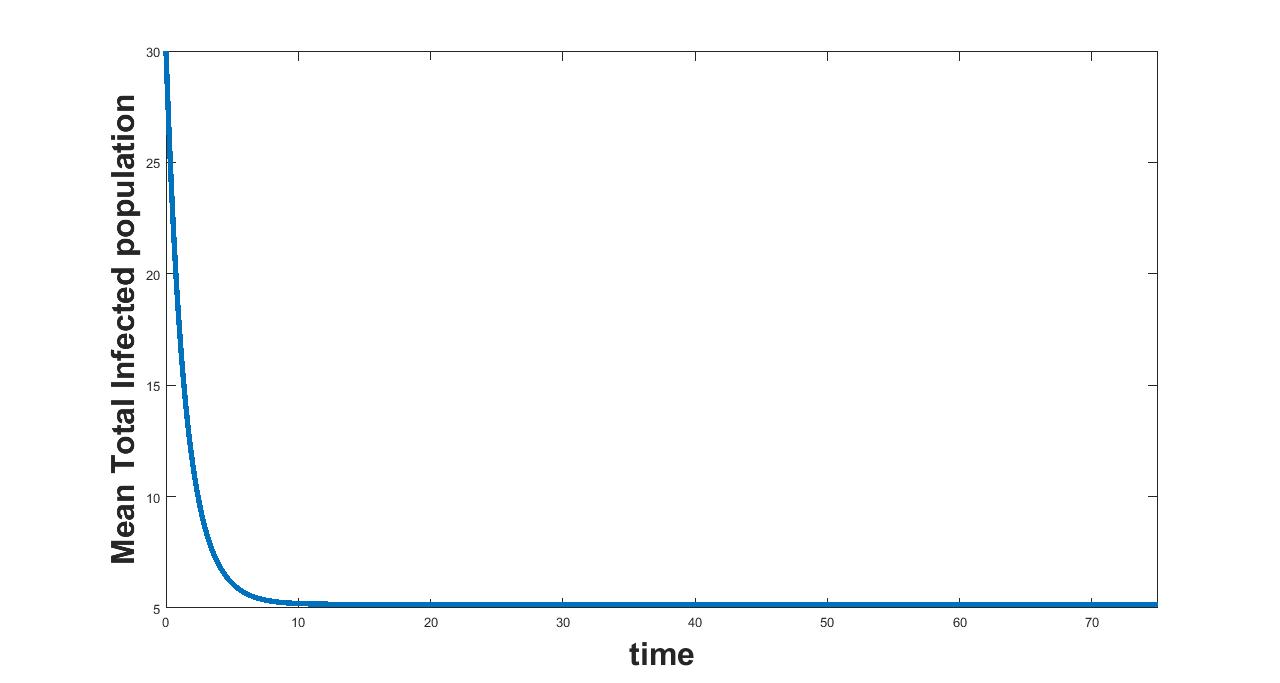}
\hspace{-.395cm}
\includegraphics[width=2.2in, height=1.8in, angle=0]{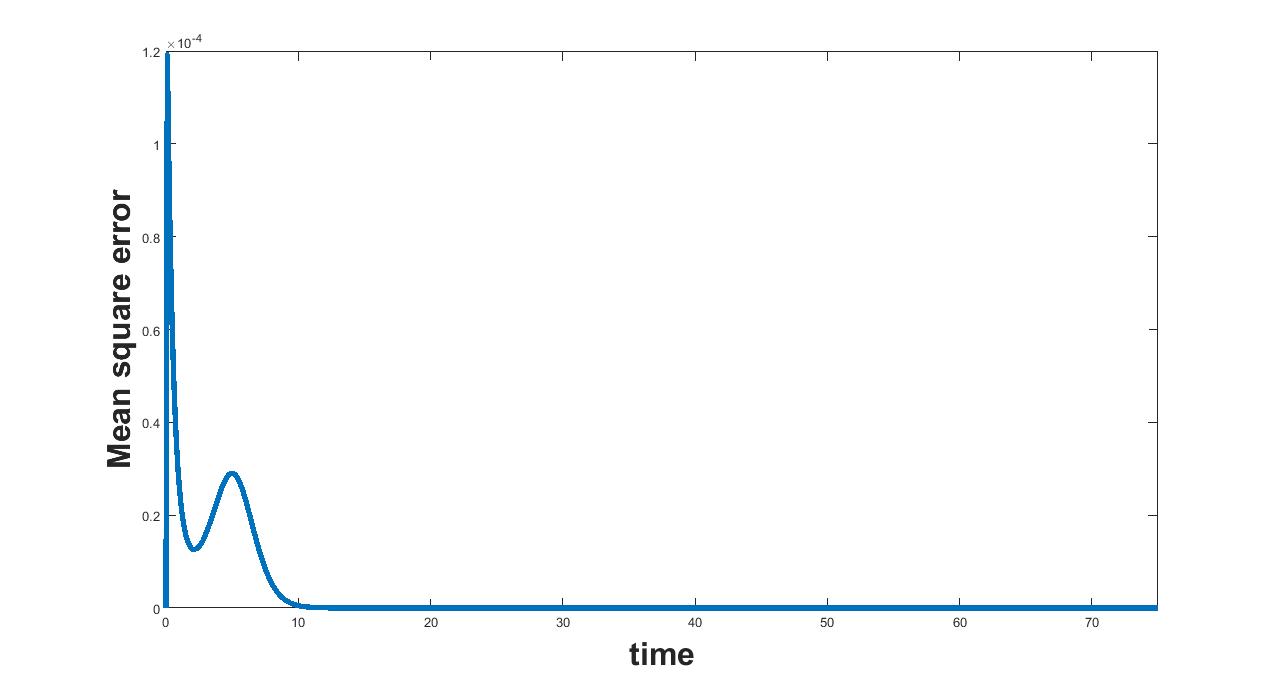}
\caption*{(b) Interval I : .0001 to .0005}
\vspace{5mm}

\end{center}
\end{figure}
\begin{figure}[hbt!]
\begin{center}
\includegraphics[width=2.2in, height=1.8in, angle=0]{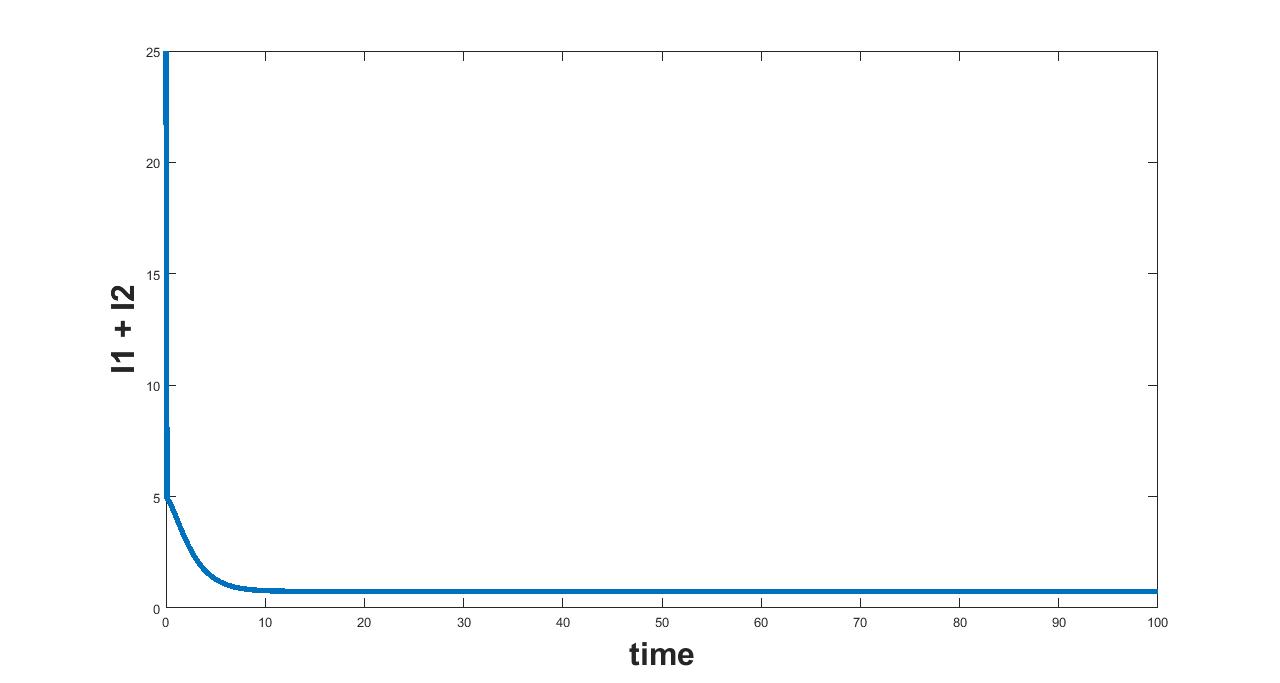}
\hspace{-.4cm}
\includegraphics[width=2.2in, height=1.8in, angle=0]{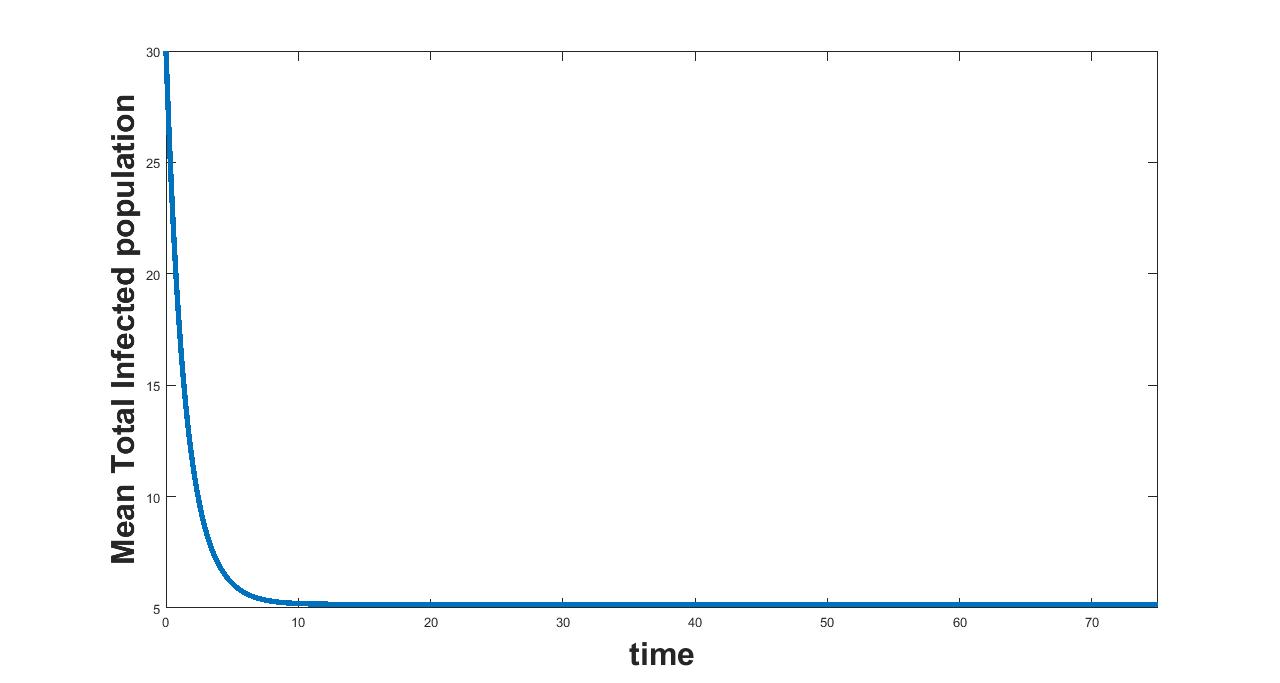}
\hspace{-.395cm}
\includegraphics[width=2.2in, height=1.8in, angle=0]{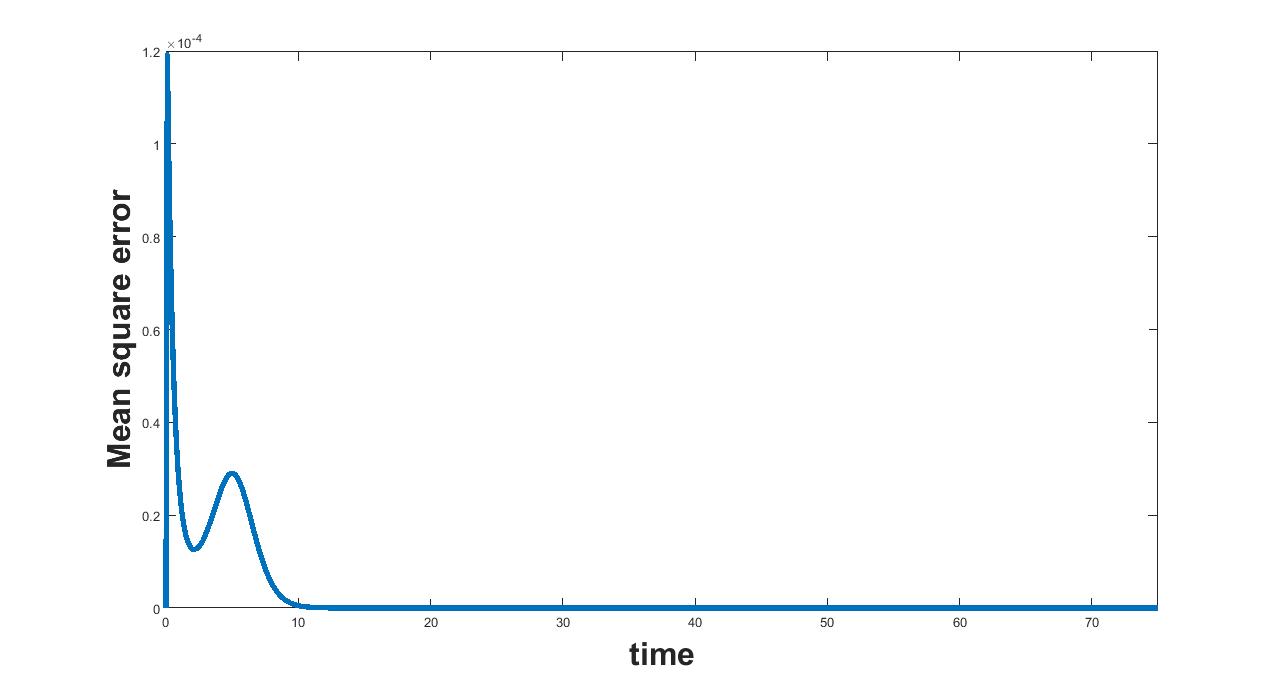}
\caption*{(b) Interval I : .0005 to .001}
\vspace{5mm}

\caption{Figure showing that the sensitivity analysis of $\beta$ varied in 2 intervals in table \ref{t}. The plots shows that  the infected population for each  value of the parameter $\beta$  per every  interval and  with  the mean infected population and the mean square error in the same interval.  }
\end{center}
\end{figure}
		
\cleardoublepage

\underline{Parameter $\boldsymbol b_1$}
The results related to sensitivity of $\boldsymbol b_1$, varied in two intervals as mentioned in table \ref {t}, are given in figure 4.12. The plots of infected population for each varied value of the parameter $\boldsymbol b_1$  per interval, the mean infected population and the mean square error are used to determine the sensitivity. We conclude from these plots that the parameter$\boldsymbol b_1$  is insensitive in interval I and  II. In similar lines, the sensitivity analysis is done for other parameters. The results are summarized below. 
\begin{figure}[hbt!]
\begin{center}
\includegraphics[width=2.2in, height=1.8in, angle=0]{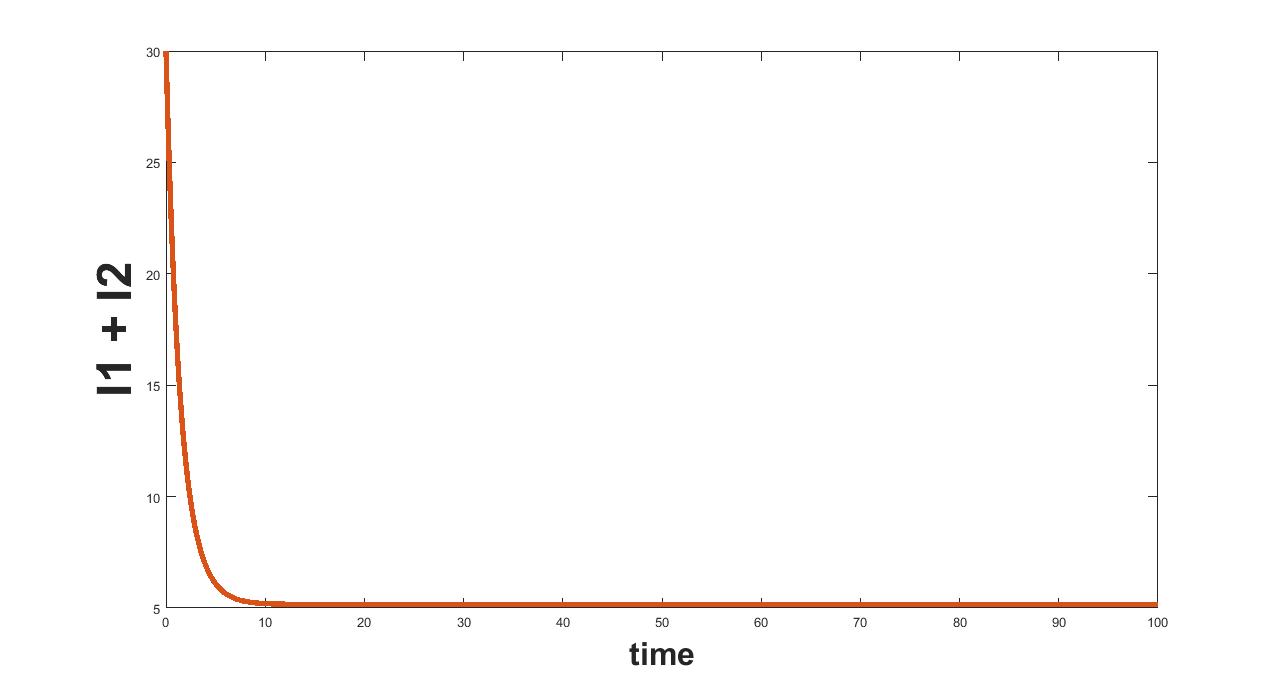}
\hspace{-.4cm}
\includegraphics[width=2.2in, height=1.8in, angle=0]{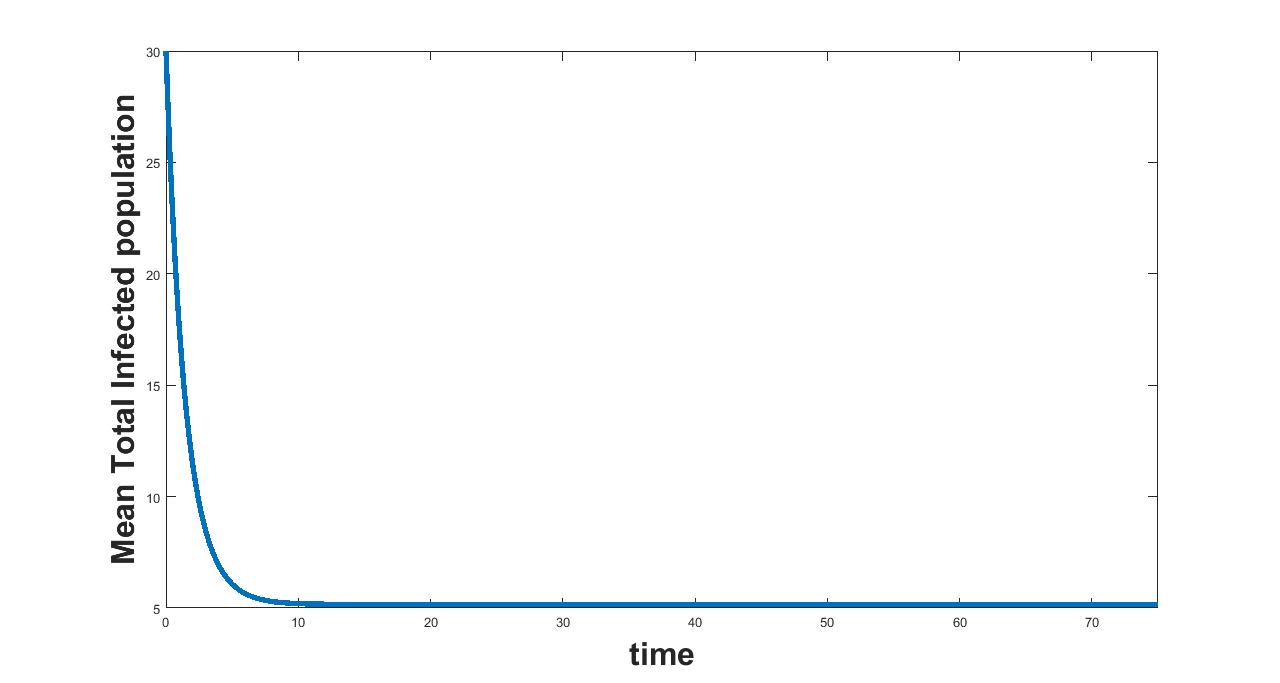}
\hspace{-.395cm}
\includegraphics[width=2.2in, height=1.8in, angle=0]{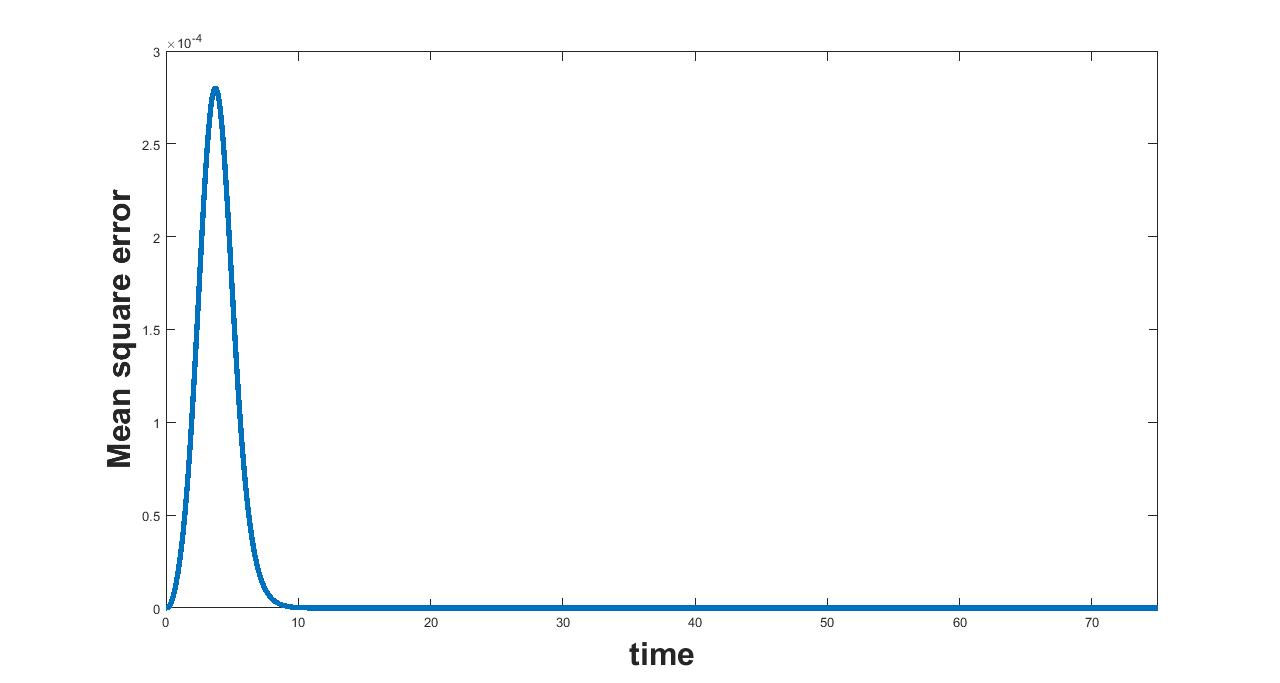}
\caption*{(a) Interval I : 345 to 355}
\end{center}
\end{figure}

\vspace{-3mm}

\begin{figure}[hbt!]
\begin{center}
\includegraphics[width=2.2in, height=1.8in, angle=0]{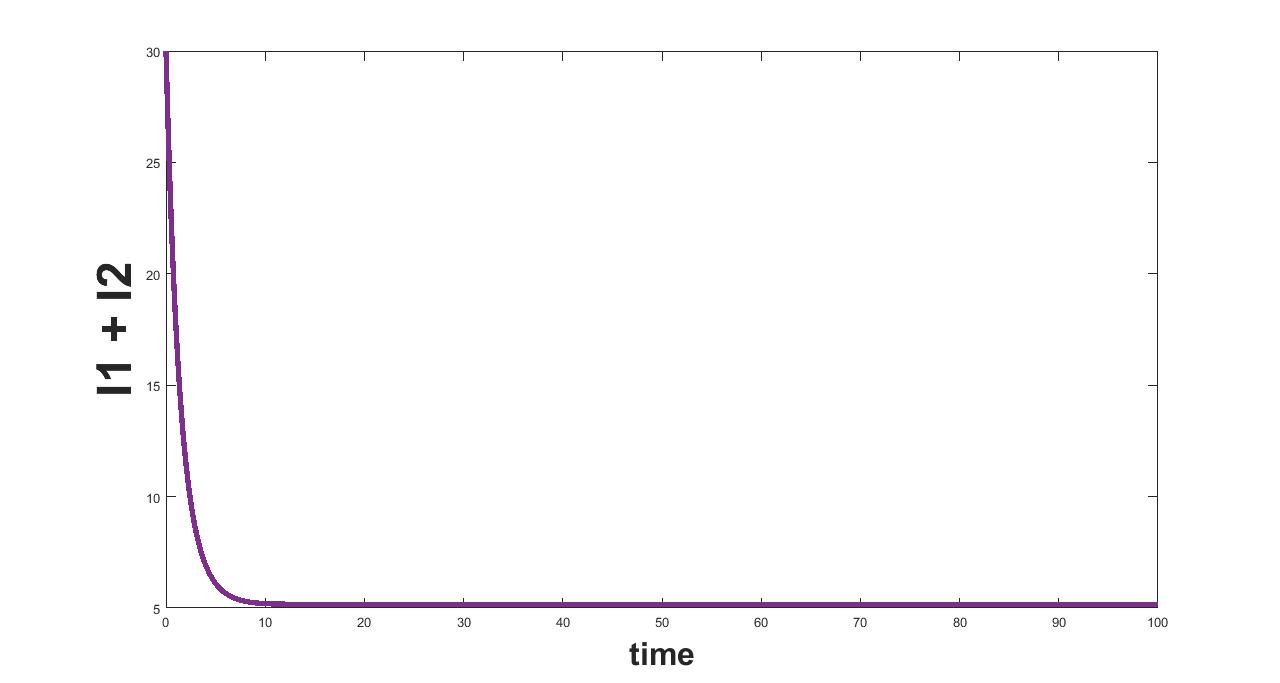}
\hspace{-.4cm}
\includegraphics[width=2.2in, height=1.8in, angle=0]{Chapters/alpha_mean.jpg}
\hspace{-.395cm}
\includegraphics[width=2.2in, height=1.8in, angle=0]{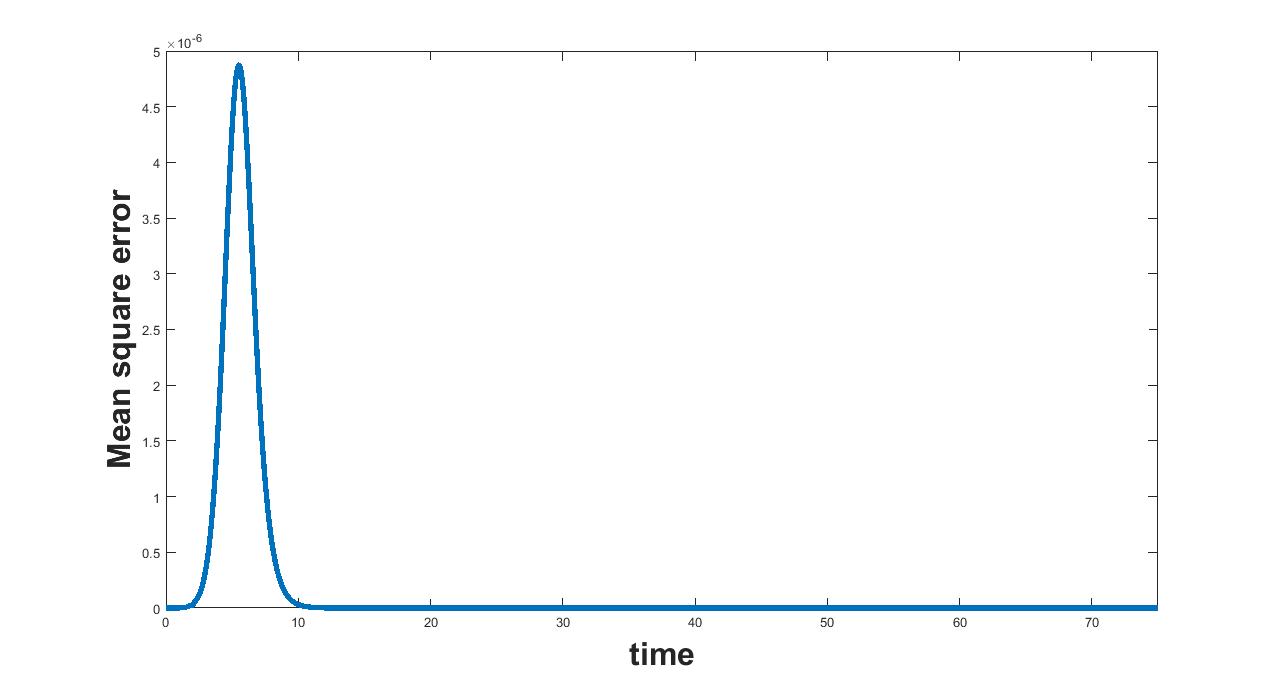}
\caption*{(b) Interval I : 355 to 365}
\vspace{5mm}

\caption{Figure showing that the sensitivity analysis of $b_1$ varied in 2 intervals in table \ref{t}. The plots shows that  the infected population for each  value of the parameter $b_1$  per every  interval and  with  the mean infected population and the mean square error in the same interval.   }
\end{center}
\end{figure}

\cleardoublepage

\underline{Parameter $\boldsymbol{m}$}
The results related to sensitivity of $\boldsymbol {m}$, varied in two intervals as mentioned in table \ref {t}, are given in figure \ref {m}. The plots of infected population for each varied value of the parameter $\boldsymbol {m}$  per interval, the mean infected population and the mean square error are used to determine the sensitivity. We conclude from these plots that the parameter$\boldsymbol {m}$  is insensitive in interval I and  II. 
	
\begin{figure}[hbt!]
\begin{center}
\includegraphics[width=2.2in, height=1.8in, angle=0]{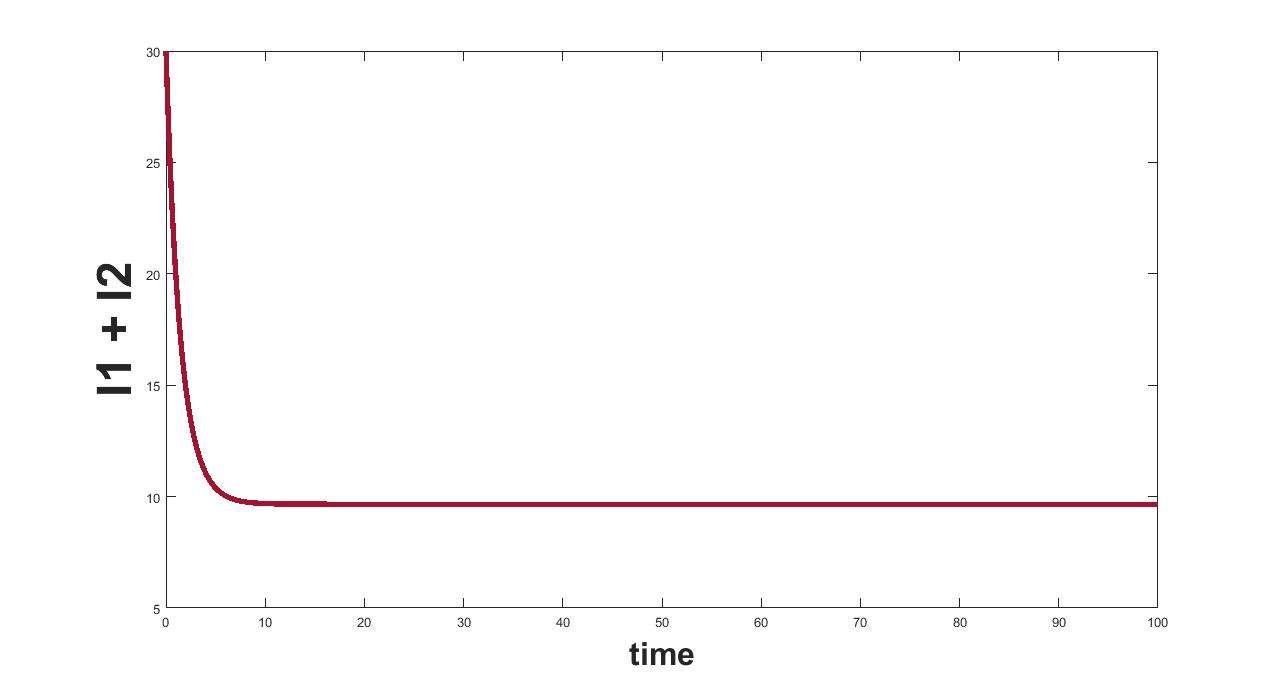}
\hspace{-.4cm}
\includegraphics[width=2.2in, height=1.8in, angle=0]{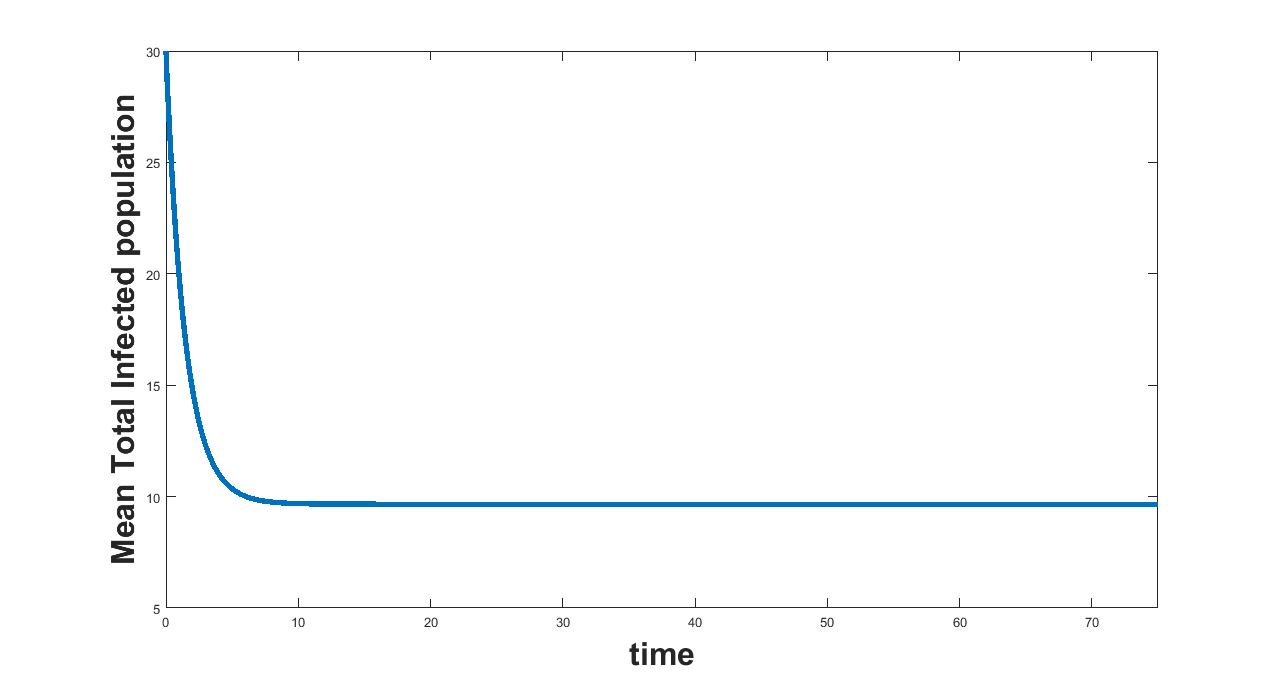}
\hspace{-.395cm}
\includegraphics[width=2.2in, height=1.8in, angle=0]{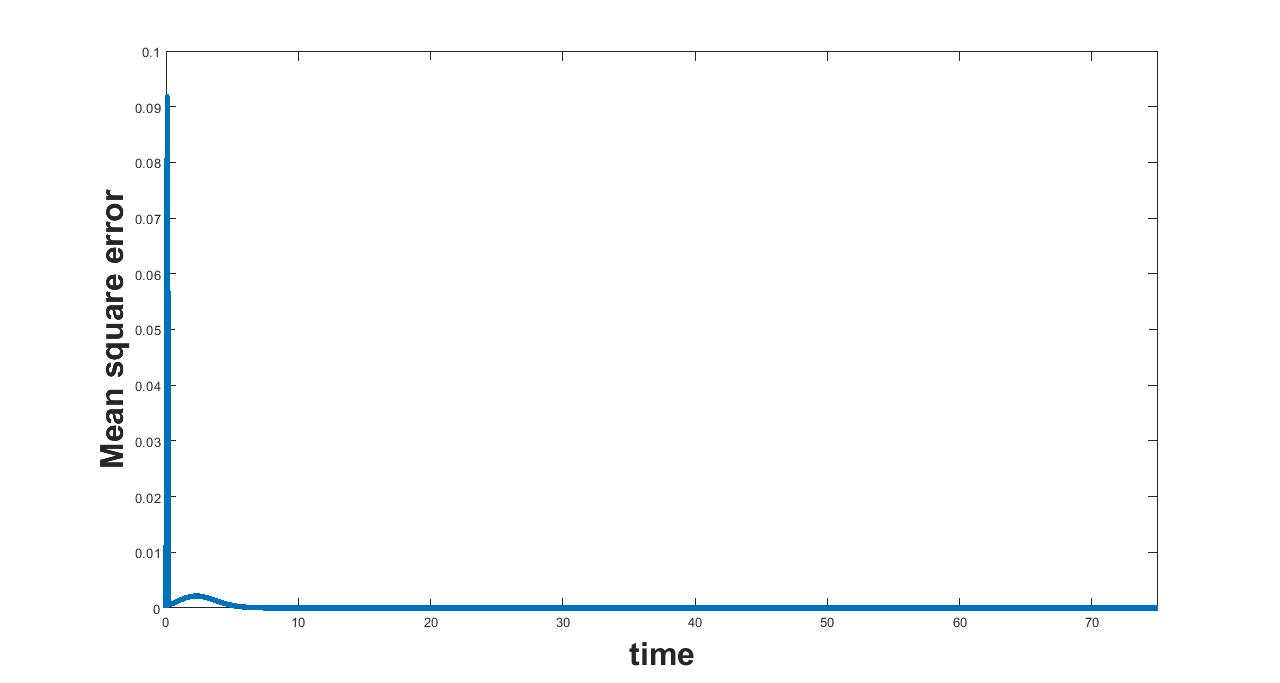}
\caption*{(a) Interval I : 0 to 2}
\end{center}
\end{figure}

\vspace{-3mm}

\begin{figure}[hbt!]
\begin{center}
\includegraphics[width=2.2in, height=1.8in, angle=0]{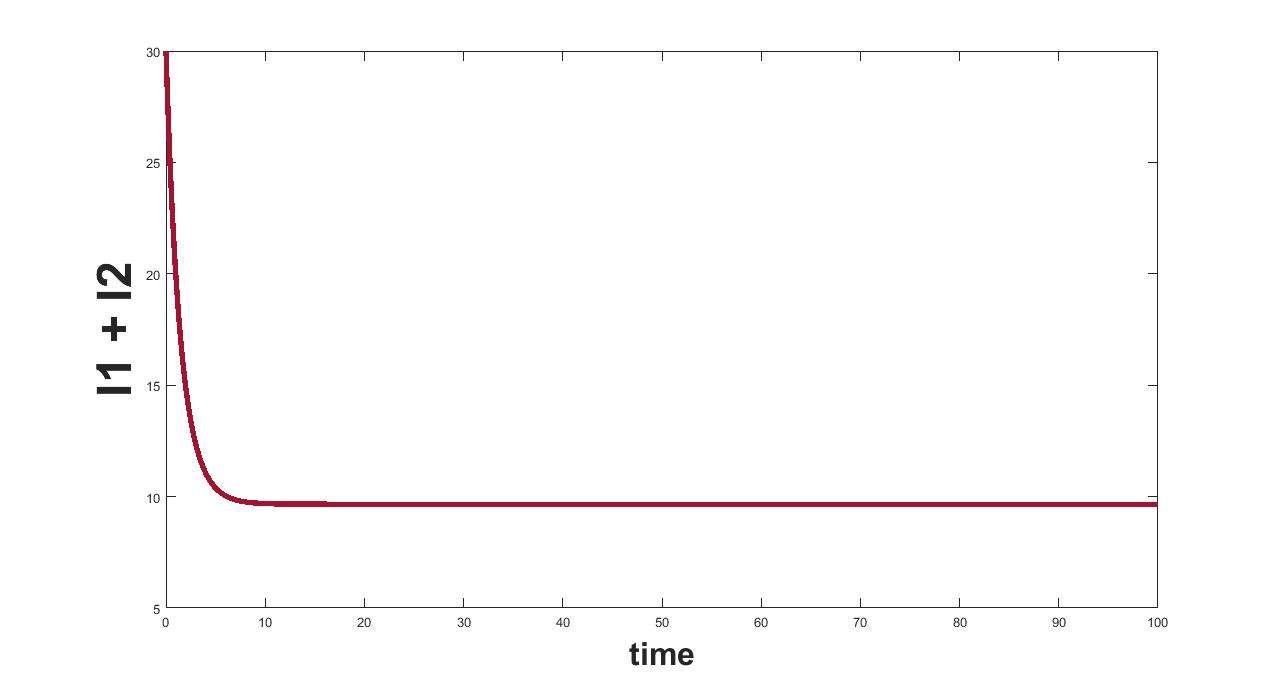}
\hspace{-.4cm}
\includegraphics[width=2.2in, height=1.8in, angle=0]{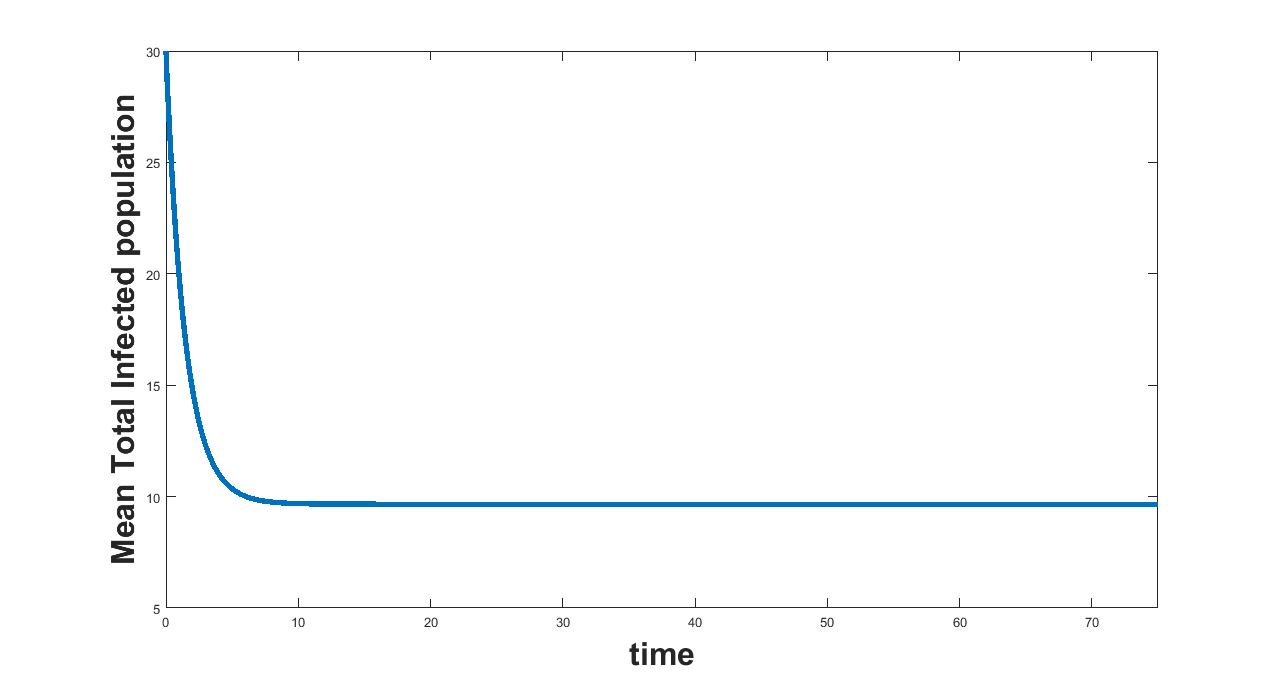}
\hspace{-.395cm}
\includegraphics[width=2.2in, height=1.8in, angle=0]{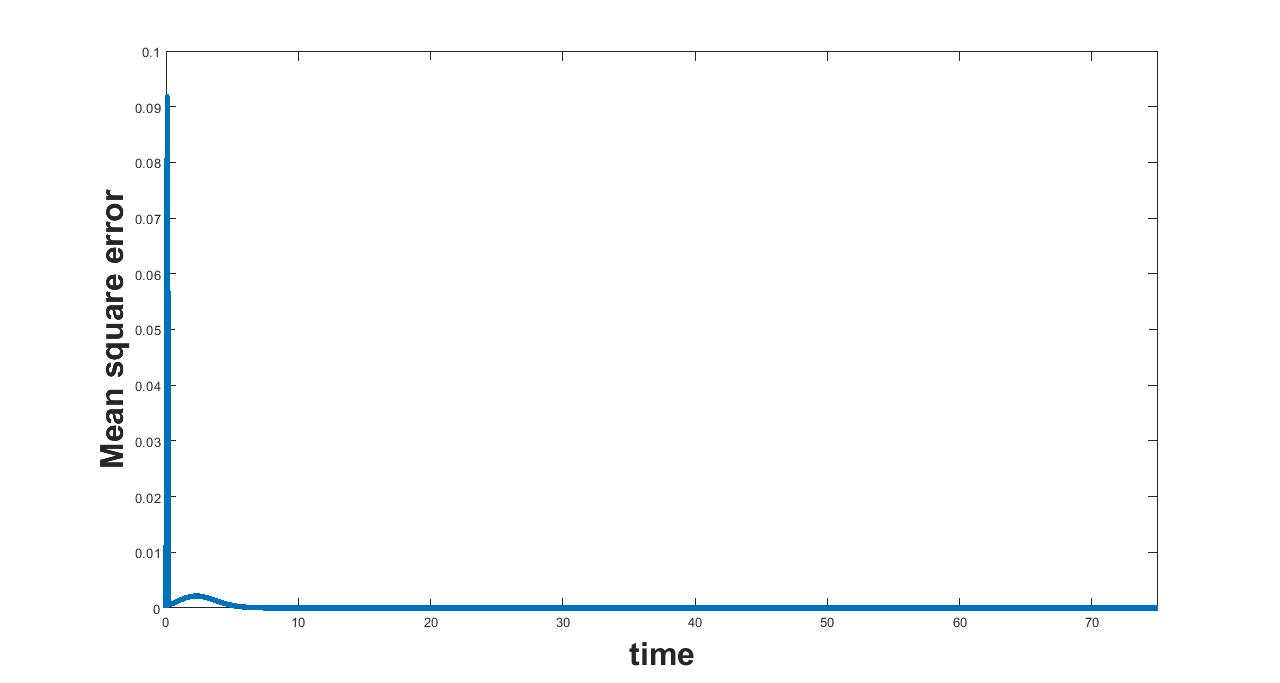}
\caption*{(b) Interval I : 2 to 3}

\vspace{5mm}
\caption{Figure showing that the sensitivity analysis of $m$ varied in 2 intervals in table \ref{t}. The plots shows that  the infected population for each  value of the parameter $m$  per every  interval and  with  the mean infected population and the mean square error in the same interval. }
\label{m}
\end{center}
\end{figure}

\cleardoublepage

 \section{Summary of Sensitivity Analysis}

The following table 4.5 gives a summary of the sensitive analysis. Parameters $\gamma$, $\mu_{c},$ $k$ are found to be sensitive in certain intervals and parameters are found to be insensitive.  \\
\begin{center}
\begin{table}[hbt!]
		\caption{Summary of Sensitivity Analysis}
		\centering
		\label{sen_anl}
	   {
			\begin{tabular}{|l|l|l|}
				\hline
				\textbf{Parameter} & \textbf{Interval} & \textbf{Step Size}  \\
				\hline

				 $b_1$ & 345 to 355  & $\times$
					\\ \cline{2-3}
				 & 355  to 365  & $\times$\\
				 	
				 \hline
				 	$ m $ & 0 to 0.00182 & $\times$
					\\ \cline{2-3}
				 & .00182 to 1  & $\times$\\
				 \hline

				$\beta$ & 0 to 0.00028 & $\times$
				\\ \cline{2-3}
				& 0.0028 to 0.1 & $\times$ \\ 
				\hline 
				
				\hline 
				$k$ & 0 to 0.05 &  \checkmark 
				\\ \cline{2-3}
				& 0.05 to 2 &  \checkmark\\
				\hline 
				$d_1$ & 0 to 0.013  & $\times$ 
				\\ \cline{2-3}
				& 0.013 to 0.5 &  $\times$ \\
				\hline 
				
				$\mu_{c}$ & 0.1 to 0.5 & \checkmark
				\\ \cline{2-3}
				& 0.5 to 1 & \checkmark \\
			 
				\hline 
			$\gamma$ & 0 to 0.0714 & \checkmark
				\\ \cline{2-3}
				& 0.0714 to 1 & \checkmark\\
				\hline 
			
			\end{tabular}
		}
\end{table} 
\end{center}

\cleardoublepage

\section{Heat Plots}

In this section we vary two sensitive  model parameters at a time in the interval given in table \ref{sen_anl} and plot the value of ${\mathcal R}_{0}$ as heat plots. The blue colour in these plots corresponds to the region where ${\mathcal R}_{0} < 1$, Therefore, for the choice of  parameters in this region, the disease free equilibrium is globally asymptotically  stable. Similarly, the green colour in these plots corresponds to the region where ${\mathcal R}_{0} > 1.$ Therefore, for the choice of  parameters in this region, the infected equilibrium is globally asymptotically stable.

\newpage
\underline{Parameters $\boldsymbol{ \mu_{c}\text { and } k}$}

\begin{figure}[hbt!]
\includegraphics[scale=0.50]{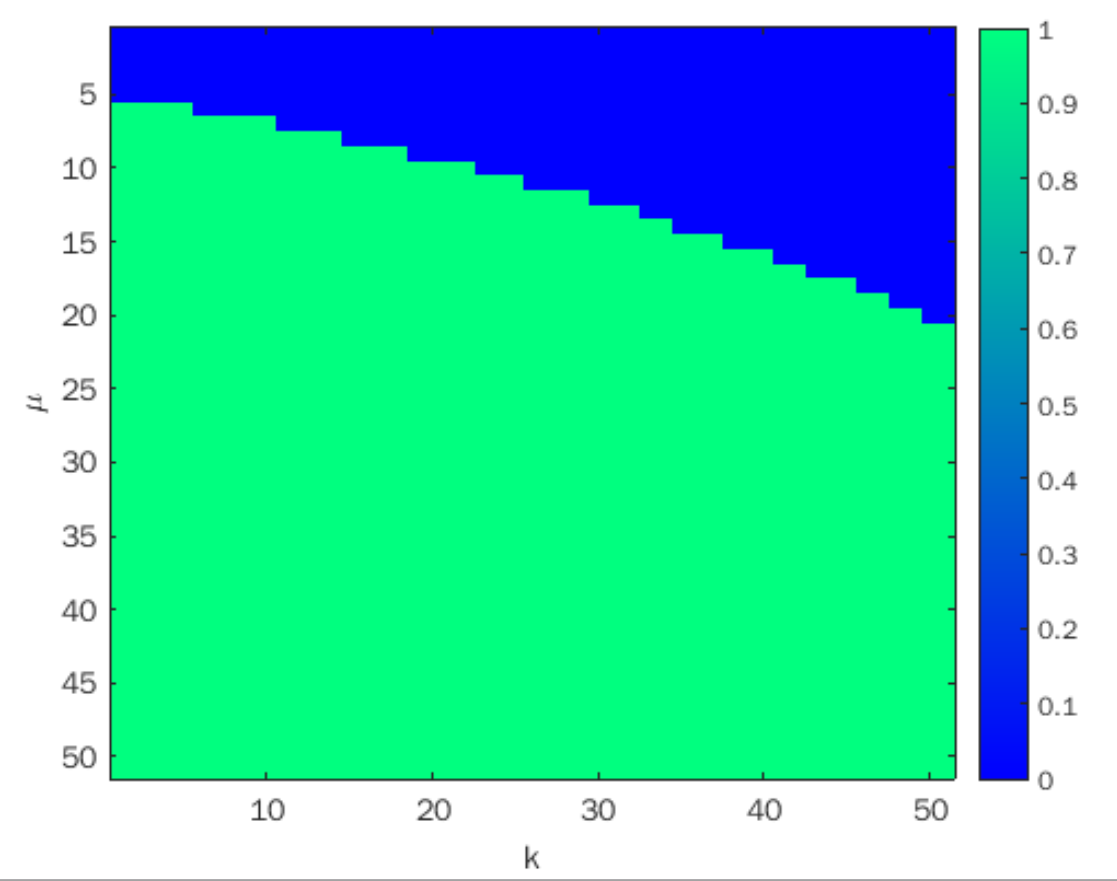}
\caption{Heat plots for the sensitive parameters $\mu_{c}$ and $k$} 
\label{hp2}
\end{figure}

\newpage

\underline{Parameters $\boldsymbol{\beta \text{ and } m}$}

\begin{figure}[hbt!]
\includegraphics[scale=0.50]{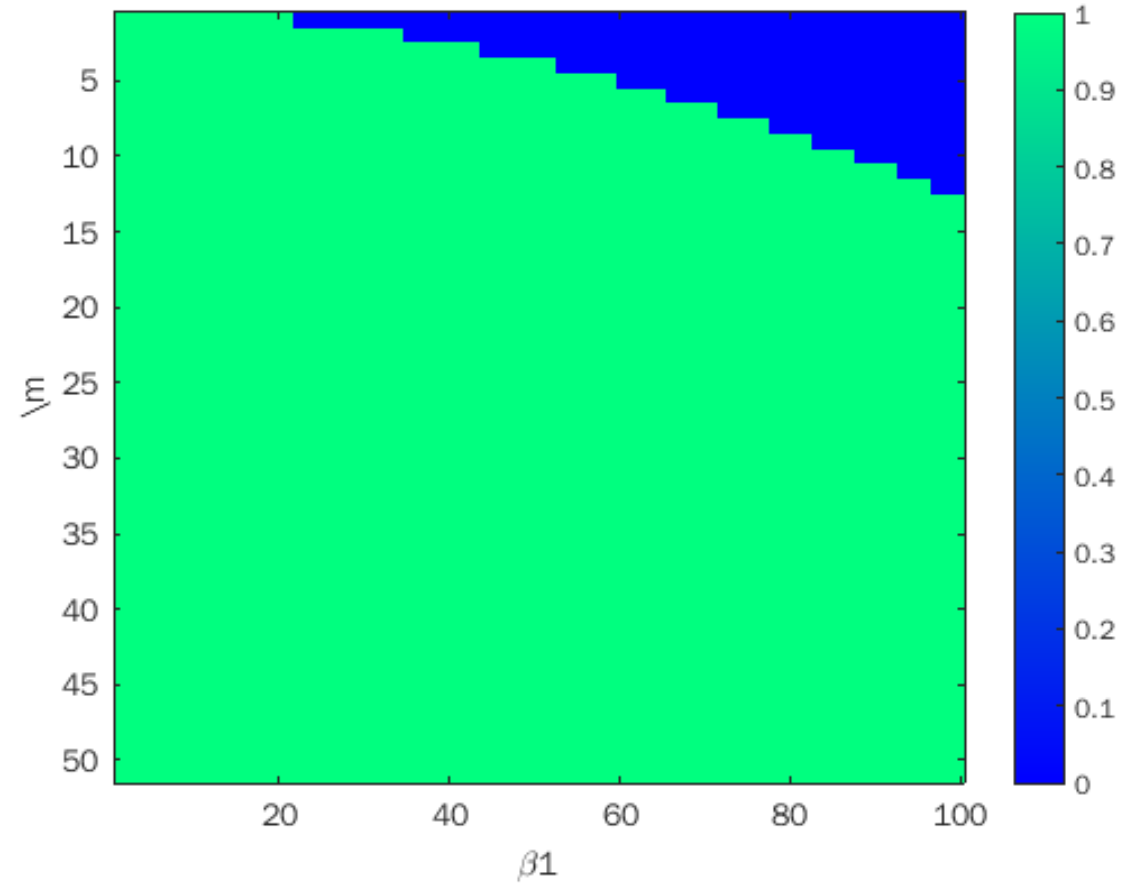}
\caption{Heat plots for the sensitive parameters $\beta$ and $m$} 
\label{hp3}
\end{figure}

\newpage
\underline{Parameters $\boldsymbol {\beta \text{ and } \mu_{c}}$}

\begin{figure}[hbt!]
\includegraphics[scale=0.50]{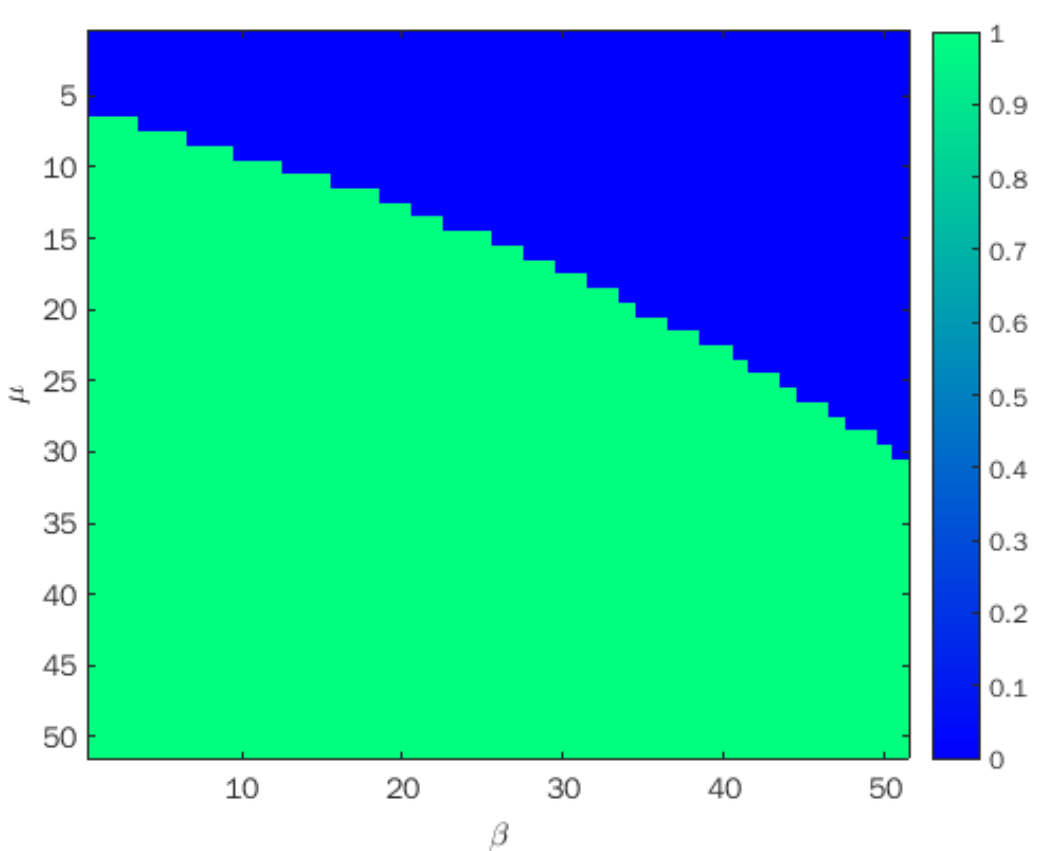}
\caption{Heat plots for the sensitive parameters $\beta$ and $\mu_{c}$} 
\label{hp4}
\end{figure}
\newpage 
\chapter [Comparative Effectiveness Study ]{\hyperlink{toc}{Comparative Effectiveness Study}}
\thispagestyle{empty}

In this chapter we do the comparative effectiveness study fr our proposed model with reference to three control interventions, namely, Vaccination, Antiviral drugs and Immunotherapy.

\newpage

{\large \bf{Vaccination}}\\

This vaccination in the population is one intervention that reduces the number of infections in both the urban and rural populations. Because this varies between urban and rural settings, we choose k to be k($1- \epsilon_{11}$) for the urban population and k($1- \epsilon_{12}$) for the rural population.\\

{\large \bf{Antiviral drugs}}\\

Antiviral medications that block viral replacement, such as Nitazoxanide, Ribavirin, and Ivermectin, aid in the reduction of COVID-19 in infected cells within seven days. Remdesivir, steroids, tocilizumab, favipiravir, and ivermectin, on the other hand, limit viral replication in infected cells \cite{av1, av2}. This gives the virus a boost, allowing it to go from infected to recovered. So  in this intervention we choose $\gamma$ to be $\gamma(1+\epsilon_{21})$ for urban population and $\gamma(1+\epsilon_{22})$ for rural population. \\

{\large \bf{ Immunotherapy}}\\

Immune-based viral elimination employing polyclonal convalescent plasma or human monoclonal antibodies to the SARS-CoV-2 spike protein may prevent infection in COVID-19-infected people or improve their outcomes \cite{imm1, imm2}. Because of antibodies, the rate of virus clearance increases as the number of infected cells decreases. As a result, the population transitions from infected to recovered (urban/rural). So here we have to choose $\gamma$ to be $\gamma(1+\epsilon_{31})$ for urban population and $\gamma(1+\epsilon_{32})$ for rural population.

    \newpage
 {\large \bf{{Change in $\mathcal{R}_{0}$}}}\\
 
\noindent
With the above  three control health interventions the modified  basic reproduction number $ \mathcal{R}_{E}$ is found to be 

$$\mathcal{R}_{E} = \frac{\beta p_1 k(1-\epsilon_{11})(1-\epsilon_{12})}{N (m+\mu_{c})(k(1-\epsilon_{11})(1-\epsilon_{12})+m+ \mu_{c})(\gamma(1+\epsilon_{21})(1+\epsilon_{22})(1+\epsilon_{31})(1+\epsilon_{32})+\mu_{c})}$$ \\

\noindent
We now do the comparative effectiveness study of these three interventions by calculating the percentage reduction of  $ \mathcal{R}_{0}$ for single and multiple combination of these interventions at different efficacy levels
such as 

\noindent
\newline(a) Low efficacy of $0.3$,
\newline (b) Medium efficacy of $0.6$, and
\newline (c) High efficacy of $0.9$. \\

The Percentage reduction of $\mathcal{R}_{0}$ is given by:

PR(percentage reduction) of $\mathcal{R}_{0}  = \bigg[ \frac{\mathcal{R}_{0} - \mathcal{R}_{E_j}}{\mathcal{R}_{0}} \bigg] \times 100$,

\noindent
where $j$ denotes for $\epsilon_{11}, \epsilon_{21}, \epsilon_{31}, \epsilon_{12}, \epsilon_{22}, \epsilon_{32}$ or the  combinations thereof seen for. \\

\noindent
For these three health interventions, we consider $17$   different efficacy combinations  as enlisted in table 5.1 for comparative effectiveness study. We did not consider other possible  efficacy combinations as they yielded more or less similar reduction in ${\mathcal R}_{0}$ as one of these 17 combinations.

\noindent
The percentages are then ranked for each efficacy level decreases on ${\mathcal R}_{0}$ for the corresponding distinct combinations for the three health interventions that were investigated in this study in ascending order from 1 to 17 for each combination of efficacy levels. CE (comparative effectiveness) is calculated and assessed on a scale of $1$to $17$, with $1$ indicating the lowest comparable effectiveness and $17$ indicating the highest comparative efficacy. The term "Comparative Effectiveness" is abbreviated as CE in table 5.2.

\newpage  

\begin{table}[htp!]
\begin{center}
\begin{tabular}{ | c | c| c | c | c | c | c | }
\hline
 \textbf{No.} & \textbf{$\epsilon_{11}$} & \textbf{$\epsilon_{12}$} & \textbf{$\epsilon_{21}$} & \textbf{$\epsilon_{22}$} & \textbf{$\epsilon_{31}$} & \textbf{$\epsilon_{32}$} \\ 
 \hline \hline
 1 &0&0&0&0&0&0    \\
 \hline
 2 & 0&0&0.3&0.3&0.3&0.3     \\
 \hline
 3 &0&0&0.3&0.6&0.3&0.6  \\ 
 \hline
 4 & 0&0&0.6&0.6&0.6&0.6   \\ 
 \hline
 5 & 0&0&0.9&0.6&0.9&0.6   \\
 \hline
 6 &0.3&0.3&0.3&0.3&0.3&0.3    \\
 \hline
 7 & 0.3&0.3&06.&0.3&0.6&0.3  \\
 \hline
 8 & 0.3&0.3&0.6&0.6&0.6&0.6    \\ 
 \hline
 9 &0.3&0.3&0.9&0.6&0.9&0.6    \\ 
 \hline
 10 & 0.6&0.6&0.3&0.3&0.3&0.3   \\ 
 \hline
 11 & 0.6&0.6&0.6&0.3&0.6&0.3  \\ 
 \hline
 12 & 0.6&0.6&0.6&0.6&0.6&0.6    \\ 
 \hline
 13 & 0.6&0.6&0.9&0.6&0.9&0.6   \\
 \hline
 14 & 0.9&0.9&0.3&0.3&0.3&0.3    \\ 
 \hline
 15 & 0.9&0.9&0.6&0.3&0.6&0.3   \\
 \hline
 16 & 0.9&0.9&0.6&0.6&0.6&0.6     \\ 
 \hline
 17 &  0.9&0.9&0.9&0.6&0.9&0.6   \\

 \hline\hline
\end{tabular}
\caption{Efficacy combinations used for CE study}
\end{center}
\end{table}

\cleardoublepage

\begin{table}[htp!]
\begin{center}
\begin{tabular}{ | c | c | c | c |}
\hline
 \textbf{No.} & \textbf{Intervention} & \textbf{\%age change in $\mathcal{R}_{0}$ } & \textbf{CE}  \\ 
 \hline \hline
 1 & $\mathcal{R}_{0}$ & 0  &  1   \\
 \hline
 2 & $\epsilon_{21}\epsilon_{22}\epsilon_{31}\epsilon_{32}$ & 4.00  & 2    \\
 \hline
 3 & $\epsilon_{21}\epsilon_{22}\epsilon_{31}\epsilon_{32}$ & 6.82  &  3  \\ 
 \hline
 4 & $\epsilon_{21}\epsilon_{22}\epsilon_{31}\epsilon_{32}$ & 10.90  &  4    \\ 
 \hline
 5 & $\epsilon_{21}\epsilon_{22}\epsilon_{31}\epsilon_{32}$ & 15.35  &  8   \\
 \hline
 6 & $\epsilon_{11}\epsilon_{12}\epsilon_{21}\epsilon_{22}\epsilon_{31}\epsilon_{32}$ & 7.70 &  4    \\
 \hline
 7 & $\epsilon_{11}\epsilon_{12}\epsilon_{21}\epsilon_{22}\epsilon_{31}\epsilon_{32}$ & 10.50  &  6    \\
 \hline
 8 & $\epsilon_{11}\epsilon_{12}\epsilon_{21}\epsilon_{22}\epsilon_{31}\epsilon_{32}$ & 14.4  &  7   \\ 
 \hline
 9 & $\epsilon_{11}\epsilon_{12}\epsilon_{21}\epsilon_{22}\epsilon_{31}\epsilon_{32}$ & 18.7  &  9   \\ 
 \hline
 10 & $\epsilon_{11}\epsilon_{12}\epsilon_{21}\epsilon_{22}\epsilon_{31}\epsilon_{32}$ & 20.35 &  10    \\ 
 \hline
 11 & $\epsilon_{11}\epsilon_{12}\epsilon_{21}\epsilon_{22}\epsilon_{31}\epsilon_{32}$ & 22.75  &  11  \\ 
 \hline
 12 & $\epsilon_{11}\epsilon_{12}\epsilon_{21}\epsilon_{22}\epsilon_{31}\epsilon_{32}$ & 26.13  &  12    \\ 
 \hline
 13 & $\epsilon_{11}\epsilon_{12}\epsilon_{21}\epsilon_{22}\epsilon_{31}\epsilon_{32}$ & 29.83   &  13  \\
 \hline
 14 & $\epsilon_{11}\epsilon_{12}\epsilon_{21}\epsilon_{22}\epsilon_{31}\epsilon_{32}$ & 80.35 &  14    \\ 
 \hline
 15 & $\epsilon_{11}\epsilon_{12}\epsilon_{21}\epsilon_{22}\epsilon_{31}\epsilon_{32}$ & 81.00  &  15    \\
 \hline
 16 & $\epsilon_{11}\epsilon_{12}\epsilon_{21}\epsilon_{22}\epsilon_{31}\epsilon_{32}$ & 81.77  &  16    \\ 
 \hline
 17 & $\epsilon_{11}\epsilon_{12}\epsilon_{21}\epsilon_{22}\epsilon_{31}\epsilon_{32}$ & 82.68  &  17    \\

 \hline\hline
\end{tabular}
\caption{Comparative Effectiveness for ${\mathcal{R}}_{0}$}
\end{center}
\end{table}

The findings of the comparative effectiveness analysis suggest the following recommendations.

\begin{itemize}
    \item [1.]  The best optimal reduction in the reproduction number was obtained when the efficacy levels of the controls intervention were chosen to be a mix of high and medium levels. 
    
     \item [2.] It is not really necessary to choose all the controls intervention at the highest efficacy levels for optimal reduction in ${\mathcal R}_0.$ 
     
        \item [3.] To achieve a fairly good reduction in  ${\mathcal R}_0$ it is necessary to choose the vaccination intervention at the highest efficacy level and the other two interventions at either medium or high efficacy levels.
    \end{itemize} 

\addtocontents{toc}{\vspace{2em}}

\newpage
\chapter{Discussions and Conclusions}

\noindent
In this study, we have formulated and analyzed a non-linear multi compartmental  (SEIR) model for the  COVID-19 with reference to immigration from urban to rural population in Indian scenario. \\

\noindent
We initially established the positivity and boudedness followed by establishing the existence and uniqueness of solutions for the proposed SEIR model. We later  calculated the equilibrium points and basic reproduction number ${\mathcal R}_0.$ \\

\noindent
We then went on to numerically establish both the local and global stabilities of the obtained disease free equilibrium and local stability of infected equilibrium. We found that the disease free equilibrium was globally asymptotically stable when ${\mathcal R}_0 < 1$ and infected equilibrium was locally asymptotically stable when ${\mathcal R}_0 > 1.$ Further performing the   sensitivity analysis we identified the sensitive parameters to be  $\gamma,$ k, $\mu_{c},$in the ranges $(0.015-0.030),(0.005-0.02),(0.003-0.005) $  respectively. 2-d two parameter heat plots with respect to sensitive parameters  were done to find the parameter regions in which the system is stable.   \\

\noindent
We finally performed the comparative  effectiveness studies  with reference to the control health interventions such as vaccination, antiviral drugs, Immunotheraphy. The findings of the comparative effectiveness analysis suggested that the best optimal reduction in the reproduction number can be  achieved  when the efficacy levels of the controls intervention are chosen to be
a mix of high and medium levels. Moreover, for achieving a fairly good reduction in reproduction number  it is necessary to choose the vaccination intervention at the highest efficacy level and the other two interventions at either medium or high efficacy levels.

\backmatter
\bibliographystyle{unsrt}
\bibliography{Bibliography.bib}

\begin{thebibliography}{10}

\bibitem{BIF}
Carlos Castillo-Chavez and Baojun Song.
\newblock Dynamical models of tuberculosis and their applications.
\newblock {\em Mathematical Biosciences and Engineering}, 1(2):361--404, 2004.

\bibitem{GLB}
Zhilan~Feng Carlos Castillo-Chavez and Wenzhang Huang.
\newblock On the computation of reproduction number and its role in global
  stability.
\newblock {\em Institute for Mathematics and Its Applications},
  125(2):229--250, 2002.

\bibitem{cooper2020sir}
Ian Cooper, Argha Mondal, and Chris~G Antonopoulos.
\newblock A sir model assumption for the spread of covid-19 in different
  communities.
\newblock {\em Chaos, Solitons \& Fractals}, 139:110057, 2020.

\bibitem{ming2020breaking}
Wai-Kit Ming, Jian Huang, and Casper~JP Zhang.
\newblock Breaking down of healthcare system: Mathematical modelling for
  controlling the novel coronavirus (2019-ncov) outbreak in wuhan, china.
\newblock {\em BioRxiv}, 2020.

\bibitem{hernandez2020host}
Esteban~A Hernandez-Vargas and Jorge~X Velasco-Hernandez.
\newblock In-host modelling of covid-19 kinetics in humans.
\newblock {\em medrxiv}, pages 2020--03, 2020.

\bibitem{kiselev2021delay}
Ilya~N Kiselev, Ilya~R Akberdin, and Fedor~A Kolpakov.
\newblock A delay differential equation approach to model the covid-19
  pandemic.
\newblock {\em medRxiv}, 2021.

\bibitem{yang2020modeling}
Cong Yang, Yali Yang, Zhiwei Li, and Lisheng Zhang.
\newblock Modeling and analysis of covid-19 based on a time delay dynamic
  model.
\newblock {\em Mathematical Biosciences and Engineering: MBE}, 18(1):154--165,
  2020.

\bibitem{aronna2020model}
M~Soledad Aronna, Roberto Guglielmi, and Lucas~M Moschen.
\newblock A model for covid-19 with isolation, quarantine and testing as
  control measures.
\newblock {\em arXiv preprint arXiv:2005.07661}, 2020.

\bibitem{dhaiban2021optimal}
Ali~Khaleel Dhaiban and Baydaa~Khalaf Jabbar.
\newblock An optimal control model of covid-19 pandemic: a comparative study of
  five countries.
\newblock {\em OPSEARCH}, pages 1--20, 2021.

\bibitem{kkdjou2020optimal}
Ramses Djidjou-Demasse, Yannis Michalakis, Marc Choisy, Micea~T Sofonea, and
  Samuel Alizon.
\newblock Optimal covid-19 epidemic control until vaccine deployment.
\newblock {\em medRxiv}, 2020.

\bibitem{libotte2020determination}
Gustavo~Barbosa Libotte, Fran~S{\'e}rgio Lobato, Gustavo~Mendes Platt, and
  Ant{\^o}nio J~Silva Neto.
\newblock Determination of an optimal control strategy for vaccine
  administration in covid-19 pandemic treatment.
\newblock {\em Computer methods and programs in biomedicine}, 196:105664, 2020.

\bibitem{ndondo2021analysis}
AM~Ndondo, SK~Kasereka, SF~Bisuta, K~Kyamakya, EFG Doungmo, and RB~M Ngoie.
\newblock Analysis, modeling and optimal control of covid-19 outbreak with
  three forms of infection in democratic republic of the congo.
\newblock {\em Results in Physics}, 24:104096, 2021.

\bibitem{bentout2021age}
Soufiane Bentout, Abdessamad Tridane, Salih Djilali, and Tarik~Mohammed
  Touaoula.
\newblock Age-structured modeling of covid-19 epidemic in the usa, uae and
  algeria.
\newblock {\em Alexandria Engineering Journal}, 60(1):401--411, 2021.

\bibitem{bubar2021model}
Kate~M Bubar, Kyle Reinholt, Stephen~M Kissler, Marc Lipsitch, Sarah Cobey,
  Yonatan~H Grad, and Daniel~B Larremore.
\newblock Model-informed covid-19 vaccine prioritization strategies by age and
  serostatus.
\newblock {\em Science}, 371(6532):916--921, 2021.

\bibitem{kumar2019role}
Anuj Kumar and Prashant~K Srivastava.
\newblock Role of optimal screening and treatment on infectious diseases
  dynamics in presence of self-protection of susceptible.
\newblock {\em Differential Equations and Dynamical Systems}, pages 1--29,
  2019.

\bibitem{mandale2021dynamics}
Roshan Mandale, Anuj Kumar, D.K.K Vamsi, and Prashant~K Srivastave.
\newblock Dynamics of an infectious disease in the presence of saturated
  medical treatment of holling type iii and self-protection.
\newblock {\em Journal of Biological Systems}, pages 1--45, 2021.

\bibitem{bishal1}
Bishal Chhetri, DKK Vamsi, and Carani~B Sanjeevi.
\newblock Optimal control studies on age structured modeling of covid-19 in
  presence of saturated medical treatment of holling type iii.
\newblock {\em Differential Equations and Dynamical Systems}, pages 1--40,
  2022.

\bibitem{diekmann2010construction}
Odo Diekmann, JAP Heesterbeek, and Michael~G Roberts.
\newblock The construction of next-generation matrices for compartmental
  epidemic models.
\newblock {\em Journal of the Royal Society Interface}, 7(47):873--885, 2010.

\bibitem{kouokam2013disease}
Etienne Kouokam, Jean-Daniel Zucker, Franklin Fondjo, and Marc Choisy.
\newblock Disease control in age structure population.
\newblock {\em International Scholarly Research Notices}, 2013, 2013.

\bibitem{samui2020mathematical}
Piu Samui, Jayanta Mondal, and Subhas Khajanchi.
\newblock A mathematical model for covid-19 transmission dynamics with a case
  study of india.
\newblock {\em Chaos, Solitons \& Fractals}, 140:110173, 2020.

\bibitem{srivastav2021modeling}
Akhil~Kumar Srivastav, Mini Ghosh, Xue-Zhi Li, and Liming Cai.
\newblock Modeling and optimal control analysis of covid-19: Case studies from
  italy and spain.
\newblock {\em Mathematical Methods in the Applied Sciences},
  44(11):9210--9223, 2021.

\bibitem{av1}
Ziyi Li, Xiaojie Wang, Donglin Cao, Ruilin Sun, Cheng Li, and Guowei Li.
\newblock Rapid review for the anti-coronavirus effect of remdesivir.
\newblock {\em Drug discoveries \& therapeutics}, 14(2):73--76, 2020.

\bibitem{av2}
CM~Chu, VCC Cheng, IFN Hung, MML Wong, KH~Chan, KS~Chan, RYT Kao, LLM Poon, CLP
  Wong, Y~Guan, et~al.
\newblock Role of lopinavir/ritonavir in the treatment of sars: initial
  virological and clinical findings.
\newblock {\em Thorax}, 59(3):252--256, 2004.

\bibitem{imm1}
Eakachai Prompetchara, Chutitorn Ketloy, and Tanapat Palaga.
\newblock Immune responses in covid-19 and potential vaccines: Lessons learned
  from sars and mers epidemic.
\newblock {\em Asian Pacific journal of allergy and immunology}, 38(1):1--9,
  2020.

\bibitem{imm2}
P~Conti, G~Ronconi, AL~Caraffa, CE~Gallenga, R~Ross, I~Frydas, and SK~Kritas.
\newblock Induction of pro-inflammatory cytokines (il-1 and il-6) and lung
  inflammation by coronavirus-19 (covi-19 or sars-cov-2): anti-inflammatory
  strategies.
\newblock {\em J Biol Regul Homeost Agents}, 34(2):327--331, 2020.

\end{thebibliography}

\end{document}